\documentclass[10pt,notitlepage]{amsart} 
\usepackage{amsmath}
\usepackage{amsthm}
\usepackage{amsfonts}
\usepackage{scalerel}
\usepackage{amssymb,amscd,epsf,verbatim,mathtools,framed}
\usepackage{mathrsfs}
\usepackage{graphicx}
\usepackage{latexsym}
\usepackage{lscape}
\usepackage[colorlinks=true]{hyperref}
\hypersetup{colorlinks, citecolor=blue, filecolor=black, linkcolor=purple, urlcolor=violet}
\usepackage{epstopdf}
\usepackage{tikz}
\usetikzlibrary{calc}
\usetikzlibrary{matrix,arrows,decorations.pathmorphing}
\usepackage{tikz-cd}
\usepackage{color}
\usepackage{multirow}  
\usepackage{scalefnt}
\usepackage{fancyhdr}
\usepackage[margin=1.32in]{geometry}
\usepackage[T1]{fontenc}

\usepackage{relsize} 
\usepackage[bbgreekl]{mathbbol} 
\DeclareSymbolFontAlphabet{\mathbb}{AMSb} %to ensure that the meaning of \mathbb does not change
\DeclareSymbolFontAlphabet{\mathbbl}{bbold} 
\newcommand{\Prism}{{\mathlarger{\mathbbl{\Delta}}}}

\newcommand{\Z}{\mathbb{Z}}

\newcommand{\N}{\mathbb{N}}

\newcommand{\Q}{\mathbb{Q}}

\renewcommand{\L}{\mathbb{L}}

\newcommand{\mD}{\mathcal{D}}
\newcommand{\mE}{\mathcal{E}}

\newcommand{\mM}{\mathcal{M}}
\newcommand{\mN}{\mathcal{N}}

\newcommand{\mO}{\mathcal{O}}

\newcommand{\mS}{\mathcal{S}}

\newcommand{\sM}{\mathscr{M}}
\newcommand{\fM}{\mathfrak{M}} 
\newcommand{\fm}{\mathfrak{m}} 

\newcommand{\tu}{\textup}
\newcommand{\cl}{\overline}
\newcommand{\ul}{\underline}

\newcommand{\ra}{\rightarrow}

\newcommand{\sq}{\widetilde}

\DeclareMathOperator{\im}{im}

\DeclareMathOperator{\Ainf}{\text{A}_{\textup{inf}}}
\DeclareMathOperator{\Binf}{\textup{B}_{\textup{inf}}}

\DeclareMathOperator{\Acrys}{\text{A}_{\textup{crys}}}

\DeclareMathOperator{\Bdrp}{\text{B}_{\textup{dR}}^+}
\DeclareMathOperator{\Bdr}{\text{B}_{\textup{dR}}}

\DeclareMathOperator{\Bstp}{\text{B}_{\textup{st}}^+}
\DeclareMathOperator{\Bst}{\text{B}_{\textup{st}}}
\DeclareMathOperator{\logcrys}{logcrys}

\DeclareMathOperator{\qsyn}{qsyn}
\DeclareMathOperator{\qrsp}{qrsp}

\DeclareMathOperator{\Fil}{Fil}

\DeclareMathOperator{\lra}{\: \longrightarrow \:} 
\DeclareMathOperator{\isom}{\;\xrightarrow{\: {}_{\sim} \:} \;}

\DeclareMathOperator{\Gal}{Gal}

\newcommand{\gr}[1]{\langle {#1} \rangle} 
\DeclareMathOperator{\spec}{Spec}
\DeclareMathOperator{\spf}{Spf}
\DeclareMathOperator{\spa}{Spa}

\DeclareMathOperator{\ett}{\textup{\'et}}
\DeclareMathOperator{\ket}{\textup{k\'et}}

\newcommand{\bi}{\begin{itemize}}
\newcommand{\ei}{\end{itemize}}
\newcommand{\bt}{\begin{theorem}}
\newcommand{\et}{\end{theorem}}
\newcommand{\bbt}{\begin{theorem*}}
\newcommand{\eet}{\end{theorem*}}
\newcommand{\bp}{\begin{proposition}}
\newcommand{\ep}{\end{proposition}}
\newcommand{\bl}{\begin{lemma}}
\newcommand{\el}{\end{lemma}}
\newcommand{\bbl}{\begin{lemma*}}
\newcommand{\eel}{\end{lemma*}}
\newcommand{\bc}{\begin{corollary}}
\newcommand{\ec}{\end{corollary}}
\newcommand{\beg}{\begin{example}}
\newcommand{\eeg}{\end{example}}
\newcommand{\br}{\begin{remark}}
\newcommand{\er}{\end{remark}}
\newcommand{\bbr}{\begin{remark*}}
\newcommand{\eer}{\end{remark*}}
\newcommand{\bd}{\begin{definition}}
\newcommand{\ed}{\end{definition}}
\newcommand{\be}{\begin{enumerate}}
\newcommand{\ee}{\end{enumerate}}
\newcommand{\bex}{\begin{exercise}}
\newcommand{\eex}{\end{exercise}}
\newcommand{\bproof}{\begin{proof}}
\newcommand{\eproof}{\end{proof}}

\theoremstyle{theorem}
\newtheorem{theorem}{Theorem}[section]

\newtheorem*{theorem*}{Theorem} 
 
\theoremstyle{definition}

\newtheorem{example}[theorem]{Example}

\newtheorem{definition}[theorem]{Definition}
\newtheorem{proposition}[theorem]{Proposition}

\newtheorem{construction}[theorem]{Construction}

\newtheorem{lemma}[theorem]{Lemma}
\newtheorem{corollary}[theorem]{Corollary}

\newtheorem{remark}[theorem]{Remark}

\usepackage{etoolbox}
\patchcmd{\section}{\scshape}{\bfseries}{}{}
\makeatletter
\renewcommand{\@secnumfont}{\bfseries}
\makeatother
\title{\large $\Z_p \tu{-lattices in semistable Galois representations}$}
\author{Zijian Yao}
\email{zijianyao@uchicago.edu}
\address{Department of Mathematics, University of Chicago}

\numberwithin{equation}{section}
\begin{document}

\begin{abstract}
We show that the category of logarithmic prismatic $F$-crystals on $(\mO_K, \varpi^\N)$ is equivalent to the category of $\Z_p$-lattices in semistable $\tu{Gal}_K$-representations. We then apply our method to describe such Galois representations using linear algebraic data via various ``logarithmic'' versions of Breuil--Kisin modules. %Along the way, we show that Galois representations of finite $E(u)$-height are potentially semistable. 
\end{abstract}

\maketitle
 \thispagestyle{empty}
%\tableofcontents
 
\vspace{-0.2in}

\section{Introduction} 
Let $K$ be a complete discretely valued extension of $\Q_p$ with ring of integers $\mO_K$ and perfect residue field $k$. Let $\varpi$ be a uniformizer in $\mO_K$. In \cite{BS_crystal}, Bhatt--Scholze provides an equivalence  between the the category of prismatic $F$-crystals on $\mO_K$ and the category of $\Z_p$-lattices in crystalline $\Gal_{K}$-representations. This equivalence has been extended to %an equivalence between certain log prismatic $F$-crystals on $(\mO_K, \varpi^\N)$ and $\Z_p$-lattices in semistable $\Gal_{K}$-representations 
the semistable case by Du--Liu \cite{DuLiu} using a different, slightly less direct approach. Our first goal is to give a more direct proof of this equivalence using the method of \cite{BS_crystal}.\footnote{In fact, we consider a variant of the absolute log prismatic site studied in \cite{DuLiu}. Our category consists of more objects which should make it easier to evaluate log $F$-crystals and study them via linear algebra data.} 

\bt \label{mainthm:equivalence}
There exists a natural equivalence of categories
\[
T_{\tu{st}}:  \tu{Vect}^{\varphi} \big((\mO_K, \varpi^\N)^{\log}_{\Prism}, \mO_\Prism  \big)  \isom  \tu{Rep}^{\tu{st}}_{\Z_p} (\tu{Gal}_K)
\]
from $F$-crystals on the absolute log prismatic site of $(\mO_K, \varpi^\N)$ to the category of finite free $\Z_p$-modules $T$ equipped with a continuous $\Gal_K$-action such that $T[1/p]$ is a semistable representation. %(we will refer to such objects simply as $\Z_p$-lattices in semistable $\Gal_K$-representations). 
\et 

We refer the reader to Definition \ref{def:absolute_log_prismatic_site} for our notion of the absolute log prismatic site. In proving the main result, we give a more geometric interpretation of the functor $T_{\tu{st}}$ (compared to the approach in \cite{DuLiu}). Our method in fact allows us to establish the following result. For the setup let $\tu{Mod}^{\phi, N_{\tu{int}}}_{/\mathfrak S}$ denote the category of finite free Breuil--Kisin modules $(\fM, \varphi)$ equipped with a monodromy operator $N: \fM/u \fM \ra \fM/u \fM$ satisfying $N \varphi = p \varphi N$ (see Definition \ref{def:BK_mod_int_monodromy} for the precise definition). We emphasize that the monodromy operator is defined integrally. Let $ \tu{Rep}^{\tu{st}, \ge 0}_{\Z_p} (\tu{Gal}_K)$ denote the category of $\Z_p$-lattices in semistable $\Gal_K$-representations with non-negative Hodge--Tate weights.  

\bt\label{mainthm:BK}
There exists a natural fully faithful embedding 
\[
D_{\mathfrak S}: \tu{Rep}^{\tu{st}, \ge 0}_{\Z_p} (\tu{Gal}_K)  \lhook\joinrel\longrightarrow \tu{Mod}^{\phi, N_{\tu{int}}}_{/\mathfrak S}.
\]
\et 
 %over $W(k)$. 

In particular, this result is a strengthening of a result of Kisin \cite[Theorem 1.3.15]{Kisin_crystal}, even after passing to the isogeny categories on both sides. In fact, we will describe the essential image of the functor $D_{\mathfrak S}$ above, using a criterion in the same spirit as in  \cite{Kisin_crystal} which describes $\tu{Rep}^{\tu{st}, \ge 0}_{\Q_p} (\tu{Gal}_K)$ inside $\tu{Mod}^{\phi, N}_{/\mathfrak S} \otimes \Q_p$. As far as we know, Theorem \ref{mainthm:BK} has not yet appeared in the literature and should shed some light on the study of semistable Galois representations. 
Moreover, by reducing mod $u$ (or, equivalently, by evaluating the log prismatic $F$-crystal on the Hyodo--Kato prism $W_{\tu{HK}}$, see Example \ref{remark:Breuil_HK}), we obtain the following integral version of Fontaine's $D_{\tu{st}}$ functor. To state the result, let $W = W(k)$ be the Witt vectors of $k$  and let $\tu{MF}^{\phi, N, \tu{int}}_{/K} $ denote the category of lattices in filtered $(\phi, N)$-modules over $K$, which consists of finite free $(\phi, N)$-stable $W$-lattices inside filtered $(\phi, N)$-modules over $K$. 
%tuples $(M, \varphi, N, \Fil^\bullet)$ where $M$ is a finite free $W$-module, $\varphi$ is $\varphi_{W}$-semilinear map $M \ra M$,  $N$ is a  $M \ra M$
\bp \label{mainprop:TongLiu}
There exists a natural faithful functor 
\[
M_{\tu{st}}: \tu{Rep}^{\tu{st}, \ge 0}_{\Z_p} (\tu{Gal}_K) \lra \tu{MF}^{\phi, N, \tu{int}}_{/K} 
\]
%from $\tu{Rep}^{\tu{st}, \ge 0}_{\Z_p} (\tu{Gal}_K)$ to the category of lattices in filtered $(\phi, N)$-modules over $K$, 
such that $M_{\tu{st}}\otimes \Q_p$ recovers Fontaines functor $D_{\tu{st}}.$
\ep 
This gives a rather simple and more conceptual construction of $M_{\tu{st}}$, which has been constructed by T. Liu \cite{TongLiu_integral_MF_lattice} using a different method. However, the functor $M_{\tu{st}}$ is not fully faithful in general as already observed in \textit{loc.cit}. Nevertheless, using a structural result on $\mathfrak S$-modules,  Liu was able to prove the following  
\bc[T. Liu {\cite[Theorem 2.3]{TongLiu_integral_MF_lattice}}] Suppose that $K/\Q_p$ is finite with ramification index $e$. 
Let $ \tu{Rep}^{\tu{st}, [0, r]}_{\Z_p} (\tu{Gal}_K)$ be the category of $\Z_p$-lattices in semistable $\Gal_K$-representations with Hodge--Tate weights in $\{0, ..., r\}$ and suppose that $er \le p-2 $. Then $M_{\tu{st}}$ induces a fully faithful embedding 
\[
M_{\tu{st}}: \tu{Rep}^{\tu{st}, [0, r]}_{\Z_p} (\tu{Gal}_K) \lhook\joinrel\longrightarrow \tu{MF}^{\phi, N, \tu{int}}_{/K}.
\]
\ec

%For the rest of the introduction we fix a compatible choice of $p$-power roots $\zeta_{p^m}$ of unity and $p$-power roots $\varpi^{1/p^m}$ of $\varpi$, and write $G_\infty = \Gal_{K_\infty}$ where $K_\infty = \cup_m K (\varpi^{1/p^m})$. 
In Section \ref{section:BK} we will describe  a category   %$\tu{Mod}^{\phi, \tau_{\tu{st}}}_{/\mathfrak S}$, $\tu{Mod}^{\phi, \tau}_{/\mathfrak S}$, and $\tu{Mod}^{\phi, G_K}_{/\mathfrak S}$, which are
$\tu{Mod}^{\phi, \tau}_{/\mathfrak S}$ of $(\phi, \tau)$-modules over $\mathfrak S$, which consists of Breuil--Kisin modules with ``decorations'' coming from the Galois side.  This category can be thought of as evaluations of log prismatic $F$-crystals on $(\mO_K, \varpi^\N)$ on the \textit{Breuil--Kisin (pre-)log prism} $(\mathfrak S, u^{\N})$ together with some additional data that captures lattices in semistable $\Gal_K$-representations, and is closely related to T. Liu's theory of $(\varphi, \hat G)$-modules. In particular, we prove the following \bt \label{mainthm:phi_tau_category}
There is a natural equivalence of categories 
\[
D^{\phi, \tau}:  \tu{Rep}^{\tu{st}, \ge 0}_{\Z_p} (\tu{Gal}_K)   \isom \tu{Mod}^{\phi,\tau}_{/\mathfrak S}.
\]
\et 
This equivalence is a variant of the results in \cite{DuLiu, Caruso, Gao_Hui} (see Remark \ref{remark:relation_to_other_phi_tau}), though our proof is different from \textit{loc.cit.} It is simpler and
uses the main equivalence (Theorem \ref{mainthm:equivalence}).  %Let us further recall that, a Galois representation $V \in \tu{Rep}_{\Q_p} (\Gal_K)$ is said to \textit{have finite $E(u)$-height} if there exists a finite free Breuil--Kisin module $(\fM, \varphi)$ such that the $G_\infty$-representation $(\fM \otimes_{\mathfrak S} W(C^\flat))^{\varphi = 1}[1/p]$ induced from the $G_\infty$-action on $W(C^\flat)$ agrees with $V|_{G_\infty}$. Using Theorem \ref{mainthm:phi_tau_category} and its proof, we deduce the following consequence. 
%\bc \label{cor:Liu_conjecture} A $\Gal_K$-representation $V$ of finite $E(u)$-height is potentially semistable.  \ec
%This result was first conjectured by T. Liu (where he phrased it as a question in \cite{Tong_Liu_a_note}). Caruso attempted to give a proof in \cite{Caruso}, but it contained a rather serious gap, as observed by Ozeki (see \cite{Gao_Hui}). More recently, Gao gave a proof of Corollary \ref{cor:Liu_conjecture} in  \cite{Gao_Hui} using locally analytic vectors and deduced a variant of Theorem \ref{mainthm:phi_tau_category} from this. Both the attempt of Caruso and the proof of Gao go through the construction of the operator $N_{\nabla}$, which Kisin uses to describe the essential image of semistable represetations in the category $\tu{Mod}^{\phi, N}_{/\mathfrak S} \otimes \Q_p$ in \cite{Kisin_crystal} (see Remark \ref{remark:Kisin_category} and Remark \ref{remark:nabla}). Our observation is that, one can in fact argue in the opposite direction and give a rather simple proof of Corollary \ref{cor:Liu_conjecture} that avoids the construction of $N_{\nabla}$. 
We refer the reader to Section \ref{ss:BK_galois} for more detail and also to Diagram (\ref{eq:the_summary_diagram}) for a concise summary of the relevant objects discussed in this paper. 

Now we briefly explain the proof of Theorem \ref{mainthm:equivalence}. Our proof largely follows the strategy of \cite{BS_crystal}, but there are some additional difficulties in the semistable case. The functor $T_{\tu{st}}$ is obtained by first considering the ``\'etale realization'' of a log prismatic $F$-crystal $\mE$ over $(\mO_K, \varpi)$. This part is similar to \cite{BS_crystal} and is essentially already established in the author's joint work with Koshikawa \cite[Theorem 7.36]{logprism}, which uses a certain form of ``log quasisyntomic descent'' developed in \textit{loc.cit.} and eventually reduces to the Riemann--Hilbert correspondence. % (\cite[Proposition 3.6]{BS_crystal}).  
To show the resulting Galois representation $T_{\tu{st}}(\mE)[1/p]$ is semistable is trickier compared to the crystalline case. One difficulty is that  Fontaine's period ring $A_{\tu{st}}$ is not $p$-complete (unlike $\Acrys$), so it does not naturally live in the log prismatic site. One is thus forced to work with slightly more ``geometric'' period rings such as $\widehat A_{\tu{st}}$ (and $\widehat B_{\tu{st}}$) studied by Kato, Breuil, and others. However, it is slightly more difficult to test admissibility directly using $\widehat B_{\tu{st}}$  (see \cite[Theorem 3.3]{Breuil}). Another (related) difficulty is that, for a semistable Galois representation, the associated filtered $(\phi, N)$-module corresponds to evaluting the $F$-crystal $\mE$ on the Hyodo--Kato prism $(W, (p), 0^{\N})$, but this (pre-log) prism does not admit a map to $A_{\tu{st}}$ or $\widehat A_{\tu{st}}$ due to the presence of log structures (see the diagram in Remark \ref{remark:diagram}). We resolve the first issue by carefully considering relevant monodromy operators and get around the second issue using an observation of Breuil (see Remark \ref{remark:Breuil_HK}), which essentially uses a certain form of Dwork's trick (in the same spirit of the Hyodo--Kato isomorphism). %Let us remark that, one can get around both issues using the notion of Breuil--Kisin--Fargues modules admitting descent by Gee--Liu (see Remark \ref{remark:Gee}), however, our treatment of the monodromy operators seems more natural and plays a role in the proof of Theorem \ref{mainthm:BK}. 
Finally, to prove that the functor $T_{\tu{st}}$ is an equivalence, we first apply a construction similar to \cite{BS_crystal} (essentially due to Kisin and Berger) to obtain a vector bundle on the locus $\spa \Prism_{\mO_C} - \{p = 0\}$, and use the result of \cite{AMMN} on the Beilinson fiber sequence to descent to $(\mO_K, \varpi^\N)^{\log}_{\Prism}$. This part is a natural extension of the argument of Bhatt--Scholze in \cite{BS_crystal}.

\br[imperfect residue field] Finally, we remark that the perfectness assumption on the residue field $k$ is essentially irrelevant in our proof of the main equivalence. In particular, one can improve Theorem \ref{mainthm:equivalence} to include imperfect residue fields with \textit {finite $p$-basis}. In this setup, however, the naive generalization of Theorem \ref{mainthm:BK} is insufficient. Instead, we need to consider Breuil--Kisin modules $(\fM, \varphi)$ equipped with a monodromy operator $N$ \textit{and} a certain quasi-nilpotent connection $\nabla$ on $\fM/u\fM$ (which vanishes if $k$ is perfect). We hope to discuss this in an upcoming joint project with H. Gao, where we expect to extend the results in this article to the imperfect residue field case. In fact, when $k$ is imperfect, Theorem \ref{mainthm:equivalence} is closely related to a family version of such an equivalence. During the preparation of this article, we learned that some important cases of this family version have been established by Du--Liu--Moon--Shimizu in their upcoming preprint, where certain cases of Theorem \ref{mainthm:equivalence} in the imperfect residue field case are used as an input. 
\er

\subsection*{Notations}
%Notation:  
We largely follow \cite{BS_crystal} regarding the definitions and notations of (bounded) prisms, prismatic site, and $F$-crystals. We follow \cite{Koshikawa, logprism} for the notion of flat morphisms of monoids, log quasisyntomic descent, log prismatic site, log prismatic cohomology and its derived variant. We use $W = W(k)$ to denote the Witt vectors of $k$, where $k$ is the residue field of $K$. Let $I \subset A$ be an ideal, then $I$-completions in this paper means derived $I$-adic completions unless otherwise specified, and is denoted by $()^{\wedge}_I$. We write $A \gr{x}^{\tu{\tiny PD}}$ for the classically $p$-completed PD polynomial ring over $A$. In the Witt vectors $W(A)$, we use $[a]$ to denote the Techmuller lift of $a \in A$. 
\subsection*{Acknowledgement}
We would like to thank Tong Liu for many enlightening discussions during the preparation of this paper. It is also a pleasure to thank Matt Emerton, Mark Kisin,  Bao Le Hung, Akhil Mathew, Peter Scholze and Teruhisa Koshikawa  for helpful conversations and suggestions. 
%\newpage 

\section{Log prismatic $F$-crystals}

In this section, we make precise our notion of the absolute log prismatic site of a pre-log ring and describe $F$-crystals on this site using log quasisyntomic covers.

\subsection{The absolute log prismatic site}
We refer the reader to \cite{Koshikawa, logprism} for the notion of $\delta_{\log}$-rings, pre-log prisms, %``log prisms'' 
and log prisms. Roughly, a $\delta_{\log}$-ring is a pre-log ring $(A, \alpha: M_A \ra A)$ equipped with a map $\delta: A \ra A$ and a map $\delta_{\log}: M_A \ra A$ satisfying some additional properties (in particular, the pair $(A, \delta)$ is required to be a $\delta$-ring). A pre-log prism is a triple $(A, I, M_A)$ where $(A, M_A)$ is a $\delta_{\log}$-ring and $(A, I)$ is a prism. We say that a pre-log prism $(A, I, M_A)$ is bounded if the underlying prism $(A, I)$ is bounded and we say that it is integral (resp. saturated) if the monoid $M_A$ is integral (resp. satuarated).  A log prism is an integral pre-log prism considered up to taking the associated log structure. 
%A ``log prism'' is a pre-log prism such that $(A, I)$ is a bounded prism and $(A, M_A)$ is a log ring. 

  %Then we have the following notion of the absolute log prismatic site of $(R, P)$.
\bd[The absolute log prismatic site] \label{def:absolute_log_prismatic_site} Let $(R, P)$ be a (classically) $p$-complete integral pre-log ring.
We write $(R, P)^{\log}_{\Prism}$ for the opposite of the category of bounded pre-log prisms $(A, I, M_A)$ with a map of pre-log rings
\[(R, P) \ra (A/I, M_A)\] and equip it with the $(p, I)$-complete homologically log flat topology introduced in \cite{logprism}. In other words, covers are maps of bounded pre-log prisms $(A, I, M_A) \ra (B, IB, M_B)$ in $(R, P)^{\log}_{\Prism}$ where $A \ra B$ is $(p, I)$-completely faithfully flat and $M_A \ra M_B$ is a flat map of monoids in the sense of \cite[Definition 4.8]{Bhatt_dR} (also see \cite[Definition 2.44, Remark 2.45]{logprism}). 
%$(\spf A/I, M_A)^a \ra (\spf R, P)^a$ of log formal  schemes such that the underlying map between the monoids over $\spf A/I$ is an isomorphism. 
\ed 

If $(A, I, M_A) \ra (B, IB, M_B)$ is a cover in $(R, P)_{\Prism}^{\log}$ and $(A, I, M_A) \ra (C, IC, M_C)$ is a map of bounded pre-log prisms, then their pushout in $(R, P)_{\Prism}^{\log}$  is given by 
\[
(D, ID, M_D) = ( B \widehat \otimes_A C, ID,  M_B \oplus_{M_A} M_C)
\]
where the completion is the derived $(p, I)$-adic completion. The underlying prism $(D, ID)$ is then a bounded prism %which is derived $(p, I)$-completely faithfully flat over $(C, IC)$ 
(see \cite[Corollary 3.12]{BS}) and the map $(C, IC, M_C) \ra (D, ID, M_D)$ is a cover in 
$(R, P)^{\log}_{\Prism}$. Therefore, $
(R, P)^{\log}_{\Prism}$ indeed forms a site. We denote by $\mO_{\Prism}$ the structure sheaf on $(R, P)_{\Prism}^{\log}$ that sends $(A, I, M_A) \mapsto A$ and by $I_{\Prism}$ the sheaf that sends $(A, I, M_A) \mapsto I$. 

\br[A variant of $(R, P)^{\log}_{\Prism}$] \label{remark:variant_site} There is a slightly more natural variant of the site $(R, P)^{\log}_{\Prism}$ considered above that takes log structures instead of pre-log structures into account. Write $(X, M_X)$ for the log formal scheme $(\spf R, P)^a$, and let ${(X, M_{X})_{\Prism, \tu{fl}}^{\log}}$ denote the opposite category of bounded log prisms $(A, I, M_{\spf A}) = (A, I, M_A)^a$ equipped with a map of log formal schemes 
\[(\spf A/I, \iota^* M_{\spf A})^a \ra (X, M_X),\]
where $\iota: \spf A/I \ra \spf A$ is the closed immersion of formal schemes. We equip ${(X, M_{X})_{\Prism, \tu{fl}}^{\log}}$ with the strict ($(p, I)$-complete) flat topology. Note that it is slightly trickier to consider homologically log flat maps for log formal schemes. For this reason, we prefer to work with $(R, P)_{\Prism}^{\log}$ in this paper. 
\er 

\br[The strict site] \label{remark:strict_sites}
Both sites $(R, P)_{\Prism}^{\log}$ and ${(X, M_X)_{\Prism, \tu{fl}}^{\log}}$ are slightly different from the ones considered in \cite{Koshikawa,DuLiu} or the saturated version in \cite{logprism}. To fix notations, let us write $(R, P)_{\Prism}^{\log,\tu{str}}$ (resp.  ${(X, M_{X})_{\Prism, \tu{fl}}^{\log,\tu{str}}}$) for the full subcategory of  $(R, P)_{\Prism}^{\log}$ (resp.  ${(X, M_{X})_{\Prism, \tu{fl}}^{\log}}$)  consisting of bounded pre-log prisms $(A, I, M_A)$  (resp. log prisms $(A, I, M_{\spf A})$) such that the map $(R, P) \ra (A/I, M_A)$ induces a strict map on the associated log formal schemes (resp. the map of log formal schemes $(\spf A/I, \iota^* M_{\spf A}) \ra (X, M_X)$ is strict). Let us equip $(R, P)_{\Prism}^{\log,\tu{str}}$ (resp.  ${(X, M_{X})_{\Prism, \tu{fl}}^{\log,\tu{str}}}$) with the $(p, I)$-complete homologically log flat topology (resp. the strict $(p, I)$-complete flat topology). The strict site ${(X, M_{X})_{\Prism, \tu{fl}}^{\log,\tu{str}}}$ is what appears in \cite{Koshikawa, DuLiu}. 
\er 

The site $(R, P)_{\Prism}^{\log}$ considered in this paper has the advantage of containing more objects (since we allow the log structures to vary) and more maps (since we allow homologically log flat maps as covers) compared to its strict variants from Remark \ref{remark:strict_sites}.  This makes it easier to evaluate $F$-crystals and allows one to use quasisyntomic descent to study them.  

% Remark that when $P = 0$, the site is still different from the usual absolute prismatic site. 

\bp \label{prop:prismatic_vector_bundles_on_coordinates} Let $ \textup{Vect} ( (R, P)^{\log}_{\Prism},  \mO_{\Prism})$ denote the category of vector bundles on the ringed topos $((R, P)^{\log}_{\Prism},  \mO_{\Prism})$. 
There is a natural equivalence 
\[
\lim_{(A, I, M_A)^a \in (R, P)_{\Prism}} \textup{Vect} (A) \isom \textup{Vect} ( (R, P)^{\log}_{\Prism},  \mO_{\Prism}).
\]
Similar statements hold when $\mO_{\Prism}$ is replaced by $(\mO_{\Prism}[1/p])^\wedge_{I_{\Prism}}$ or $(\mO_{\Prism}[1/I_{\Prism}])^\wedge_{p}$. 
\ep

\bproof 
The same proof of \cite[Proposition 2.7]{BS_crystal} still applies in our setup, using the result of Drinfeld and Mathew \cite[Theorem 2.2]{BS_crystal} for the sheaves $(\mO_{\Prism}[1/p])^\wedge_{I_{\Prism}}$ and $(\mO_{\Prism}[1/I_{\Prism}])^\wedge_{p}$. 
\eproof

\subsection{Comparison to the log quasisyntomic site} \label{ss:qsyn_site}

\bd Suppose that $(R, P)$ is log quasisyntomic. 
We let ${(R, P)}_{\qsyn}$ be the opposite of the category consisting of log quasisyntomic maps $(R, P) 
\ra (S, M)$ defined in \cite{logprism}, equipped with the log quasisyntomic topology. Let 
\[{(R, P)}_{\qrsp} \subset {(R, P)}_{\qsyn}\]
be the full subcategory generated by those $(S, M)$ that are log quasiregular semiperfectoid over $(R, P)$. 
\ed 

By (the proof of) \cite[Corollary 3.20]{logprism}, we have an equivalence of categories 
\[
\textup{Shv}( {(R, P)}_{\qsyn}) \isom \textup{Shv}( {(R, P)}_{\qrsp}).
\]
For each $(S, M) \in (R, P)_{\qrsp}$, there is an ``absolute derived log prismatic cohomology'' $\Prism_{(S, M)}$ constructed in \cite[Section 4.3]{logprism}.  
Following the notation from \cite{BS_crystal}, we have the following sheaves on $(R, P)_{\qrsp}$ (thus on $(R, P)_{\qsyn}$). 
\bi
\item The structure sheaf $\Prism_{\bullet}$, which sends
\[ (S, M) \longmapsto \Prism_{(S, M)}. \]
\item The (Frobenius twisted) $\Acrys$-sheaves $\Prism_{\bullet}\{\varphi^n(I)/p\}$ for each $n \in \N$, which sends
\[ (S, M) \longmapsto \Prism_{(S, M)} \{ \varphi^n(I)/p\}. \] Note that when $n = 0$, we have $\Prism_{\bullet} \{I/p\} (S, M) \cong \Prism_{(S/p, M)}.$
\item The (Frobenius twisted) rational localization sheaves $\Prism_{\bullet} \gr{\varphi^n(I)/p} $, which sends 
\[ (S, M) \longmapsto \Prism_{(S, M)} [\varphi^n(I)/p]^{\wedge}_p. \]
\item The \'etale structure sheaf $\Prism_{\bullet} [1/I]^{\wedge}_{p}$,  which sends
\[ (S, M) \longmapsto \Prism_{(S, M)} [1/I]^{\wedge}_p. \]
\item The $\Bdrp$ (resp. $\Bdr$) period sheaf $\mathbb{B}^+_{\tu{dR}}$ (resp. $\mathbb{B}_{\tu{dR}}$), which sends  
\[(S, M) \longmapsto  \Prism_{\bullet} [1/p]^{\wedge}_{I} \quad (\tu{resp. }  \Prism_{\bullet} [1/p]^{\wedge}_{I} [1/I] )\]
\ei
All but the $\Bdrp$ and $\Bdr$ period sheaves above carry compatible Frobenius operators, which we denote by $\varphi.$ Note that, by (the proof of) \cite[Lemma 6.7]{BS_crystal}, we have natural isomorphisms 
\begin{equation} \label{eq:identifying_bdr_sheaves}
    \mathbb{B}^+_{\tu{dR}} =  \Prism_{\bullet} [1/p]^{\wedge}_{I} \isom  \Prism_{\bullet} \gr{\varphi^n(I)/p} [1/p]^{\wedge}_{I}
\end{equation}
of sheaves compatible with the restriction maps $ \Prism_{\bullet}\gr{\varphi^{n}(I)/p} \ra \Prism_{\bullet}\gr{\varphi^{n-1}(I)/p}.$

\bp \label{prop:prismatic_vector_bundle_on_qrsp}  
Let $\textup{Vect} ( (R, P)_{\qsyn},  \Prism_{\bullet})$ denote the category of vector bundles on the ringed topos $((R, P)_{\qsyn},  \Prism_{\bullet})$ (and similarly define $\textup{Vect} ( (R, P)_{\qrsp},  \Prism_{\bullet})$).  
There are natural equivalences of categories 
\[
\lim_{(S, M) \in (R, P)_{\qrsp}} \textup{Vect} (\Prism_{(S, M)}) 
\cong \textup{Vect} ((R, P)_{\qsyn}, \Prism_{\bullet})
\cong  \textup{Vect} ( (R, P)^{\log}_{\Prism},  \mO_{\Prism}).
\] 
Similar statements hold when $\mO_{\Prism}$ is replaced by $(\mO_{\Prism}[1/p])^\wedge_{I_{\Prism}}$ (resp. by $(\mO_{\Prism}[1/I_{\Prism}])^\wedge_{p}$) and   $\Prism_{\bullet}$ is replaced by $\mathbb{B}^+_{\tu{dR}}$ (resp. by $\Prism_{\bullet} [1/I]^{\wedge}_p$).
\ep 

\bproof 
The first equivalence is proved similarly as Proposition \ref{prop:prismatic_vector_bundles_on_coordinates}. The second equivalence follows from the same proof of \cite[Proposition 2.14]{BS_crystal}. Let us recall the salient point of the argument for the sheaf $\mO_{\Prism}$ (the arguments for the other sheaves are similar). For any quasisyntomic cover $(R, P) \ra (S^0, M^0)$ where $(S^0, M^0) \in (R, P)_{\tu{qrsp}}$, we can form the Cech nerve $(S^{\bullet}, M^{\bullet})$ of this cover, where each $(S^i, M^i)$ is quasiregular semiperfectoid. The absolute log prismatic cohomology $\Prism_{(S^0, M^0)}$ gives a cover of the final object in $\tu{Shv}((R, P)^{\log}_{\Prism})$ and the Cech nerve of this cover is given by $\Prism_{(S^\bullet, M^\bullet)}$, therefore the limit
\begin{equation} \label{eq:computing_crystal_unsaturated}
\lim \tu{Vect} (\Prism_{(S^\bullet, M^\bullet)})
\end{equation}
computes both $\textup{Vect} ((R, P)_{\qsyn}, \Prism_{\bullet})$
and $\textup{Vect} ( (R, P)^{\log}_{\Prism},  \mO_{\Prism})$. Taking the limit over all such covers $(R, P) \ra (S^0, M^0)$, we arrive at the desired canonical equivalence as claimed in the proposition. 
\eproof

\subsection{Log prismatic $F$-crystals}

\bd 
Let $(R, P)$ be a $p$-adically complete integral pre-log ring. An $F$-crystal on $(R, P)^{\log}_{\Prism}$ is a vector bundle $\mE \in \text{Vect} ((R, P)_{\Prism}^{\log}, \mO_{\Prism})$ equipped with an isomorphism 
\[
\varphi_{\mE}: \varphi^* \mE [1/I_{\Prism}] \isom   \mE [1/I_{\Prism}]. 
\] 
We often refer to $\mE$ as a \textit{log prismatic $F$-crystal on $(R, P)$}. 
\ed 

We denote the category of log prismatic $F$-crystals on $(R, P)$ by $ \text{Vect}^\varphi ((R, P)_{\Prism}^{\log}, \mO_{\Prism})$. More generally, for a sheaf $\mO'$ of $\mO_{\Prism}$-algebras equipped with a compatible Frobenius, we make a similar definition of the category $ \text{Vect}^\varphi ((R, P)_{\Prism}^{\log}, \mO')$ of log prismatic $F$-crystals over $\mO'$.  The category of $F$-crystals (of vector bundles) on the log quasisyntomic site $(R,P)_{\qsyn}$ (or its variant $(R, P)_{\qrsp}$) is defined similarly and denoted by   $\textup{Vect}^{\varphi}( (R, P)_{\qsyn},  \Prism_{\bullet})$ (resp. $\textup{Vect}^{\varphi}( (R, P)_{\qrsp},  \Prism_{\bullet})$). 
From   Proposition \ref{prop:prismatic_vector_bundles_on_coordinates} and  \ref{prop:prismatic_vector_bundle_on_qrsp} we obtain the following equivalences:
\begin{multline}
 \quad  \textup{Vect}^{\varphi} ( (R, P)^{\log}_{\Prism},  \mO_{\Prism}) \:\cong\:     \lim_{(A, I, M_A) \in (R, P)^{\log}_{\Prism}} \textup{Vect}^{\varphi} (A)  \\  \; \cong \; 
\lim_{(S, M) \in (R, P)_{\qrsp}} \textup{Vect}^{\varphi}(\Prism_{(S, M)}) 
\; \cong \; \textup{Vect}^{\varphi}((R, P)_{\qrsp}, \Prism_{\bullet}) .  \quad   
\end{multline}

\subsection{The \'etale realization}

Let $(R, P)$ be a $p$-adically complete pre-log ring and assume that $P$ is an fs monoid. Let $(X, M_X) = (\spf R, P)^a$ be the associated log formal scheme over $\spf \Z_p $ as in Remark \ref{remark:variant_site}. Let $X_{\eta}$ denote the generic fiber of $(X, M_X)$ over $\spa \Q_p$, viewed as a log adic space\footnote{or as a log diamond in the sense of \cite{logprism}. In fact, for the statement of Proposition \ref{prop:etale_realization}, one may simply consider the fs log scheme $(\spec R[\frac{1}{p}], P)^a$.}, and let $X_{\eta, \ket}$ denote the Kummer \'etale site of $X_{\eta}$.

\bp[{\cite[Theorem 7.36]{logprism}}] \label{prop:etale_realization}
There is a natural equivalence 
\begin{equation}\label{eq:etale_realization}
\textup{Vect}^{\varphi} ((R, P)^{\log}_{\Prism}, \mO_{\Prism}[1/I_{\Prism}]^{\wedge}_p) \isom \tu{Loc}_{\Z_p} (X_{\eta, \ket})
\end{equation}
between Laurent $F$-crystals on $(R, P)$ and $\Z_p$-valued Kummer \'etale local systems on $X_{\eta}$. 
\ep 

\bproof 
This is a slight variant of \cite[Theorem 7.36]{logprism} (see also \cite[Corollary 3.8]{BS_crystal}). For the proof, it is convenient to consider an intermediate site $(R, P)^{\log, \tu{sat}}_{\Prism}$, which is a full subcategory of $(R, P)^{\log, \tu{sat}}$ that consists of pre-log prisms $(A, I, M_A)$ where $M_A$ is saturated and the topology is the $p$-complete homologically log flat topology. Then by the same proof of \cite[Theorem 7.36]{logprism}, we have an equivalence 
\[
\textup{Vect}^{\varphi} ((R, P)^{\log, \tu{sat}}_{\Prism}, \mO_{\Prism}[1/I_{\Prism}]^{\wedge}_p) \isom \tu{Loc}_{\Z_p} (X_{\eta, \ket}). 
\]
It remains to show that restricting to saturated pre-log prisms induces a natural equivalence 
\begin{equation} \label{eq:equiv_unsat_sat}
\textup{Vect}^{\varphi} ((R, P)^{\log}_{\Prism}, \mO_{\Prism}[1/I_{\Prism}]^{\wedge}_p)  \isom \tu{Vect}^{\varphi} ((R, P)^{\log, \tu{sat}}_{\Prism}, \mO_{\Prism}[1/I_{\Prism}]^{\wedge}_p).   
\end{equation}
To this end, we may consider quasisyntomic covers $(R, P) \ra (S, M)$ of the form $(S, M) = (R_\infty, P_\infty)$ described in \cite[Section 6.4]{logprism}, where $P_\infty$ is the pushout of $\N^{\oplus J} \ra \N[1/p]^{\oplus J}$ along a surjective map $\N^{\oplus J} \ra P$ for some finite index set $J$. Let $(R, P) \ra (S^\bullet, M^\bullet)$ be the Cech nerve of such a cover, then $\textup{Vect}^{\varphi} ((R, P)^{\log}_{\Prism}, \mO_{\Prism}[1/I_{\Prism}]^{\wedge}_p)$ is computed by  
\begin{equation} \label{eq:computing_laurent_crystal_unsaturated} 
\lim \tu{Vect}^{\varphi} (\Prism_{(S^\bullet, M^\bullet)}[1/I]^{\wedge}_p)
\end{equation}
by Proposition \ref{prop:prismatic_vector_bundle_on_qrsp} and its proof. Now let us observe that, by \cite[Proposition 3.5]{BS_crystal} and the proof of \cite[Corollary 3.7]{BS_crystal}, we have natural equivalences 
\begin{align} 
\nonumber %\label{eq:perf_unsat}
\lim \tu{Vect}^{\varphi} (\Prism_{(S^\bullet, M^\bullet)}[1/I]^{\wedge}_p) &  \isom 
\lim \tu{Vect}^{\varphi} (\Prism_{(S^\bullet, M^\bullet), \tu{perf}}[1/I]^{\wedge}_p)  \\
& \label{eq:perf_p_sat} \isom 
\lim \tu{Vect}^{\varphi} (\Prism_{(S^{\bullet, \tu{p-sat}}, M^{\bullet, \tu{p-sat}}), \tu{perf}}[1/I]^{\wedge}_p)  \\
& \nonumber %\label{eq:perf_p_sat2} 
\isom 
\lim \tu{Vect}^{\varphi} (\Prism_{(S^{\bullet, \tu{p-sat}}, M^{\bullet, \tu{p-sat}})}[1/I]^{\wedge}_p),   
\end{align}
where the supscript $\tu{p-sat}$ denotes the $p$-saturation of monoids (and pre-log rings),
%For the equivalences (\ref{eq:perf_unsat}) and (\ref{eq:perf_p_sat2}) above,  
and we use \cite[Corollary 4.21]{logprism} for the equivalence in (\ref{eq:perf_p_sat}) (see also \cite[Remark 6.20]{logprism}). For each $i$, the derived log prismatic cohomology $\Prism_{(S^{i, \tu{p-sat}}, M^{i, \tu{p-sat}})}$ comes equipped with the structure of an integral pre-log prism 
\begin{equation} \label{eq:p_saturated_prism_cover}
(\Prism_{(S^{i, \tu{p-sat}}, M^{i, \tu{p-sat}})}, I, \sq M^{i, \tu{p-sat}})
\end{equation}
where $\sq M^{i, \tu{p-sat}}$ denotes the exactification of the map $(\sq M^{i, \tu{p-sat}})^{\flat} \ra \sq M^{i, \tu{p-sat}}$ (see \cite[Section 4.3, in particular, Proposition 4.14 and Proposition 4.17]{logprism}). Now we observe that $\sq M^{i, \tu{p-sat}} $ is in fact already saturated by construction (see \cite[Lemma 6.18]{logprism}), thus the homotopy limit in (\ref{eq:computing_laurent_crystal_unsaturated}) also computes $\tu{Vect}^{\varphi} ((R, P)^{\log, \tu{sat}}_{\Prism}, \mO_{\Prism}[1/I_{\Prism}]^{\wedge}_p)$. Finally, taking the limit over all such covers $(R, P) \ra (S, M)$ of the form described above, we get the desired equivalence in (\ref{eq:equiv_unsat_sat}). 
\eproof

%\subsection{Relation to the log crystalline site}

%\section{Log prismatic $F$-crystals on $(\mO_K, \varpi^{\N})$} 
\section{$F$-crystals on $(\mO_K, \varpi^{\N})^{\log}_{{\large \mathbbl{\Delta}}}$} \label{section:F_crystals}

In this section we specialize our discussion to $(\mO_K, \varpi^{\N})$. For simplicity let us write $N_\infty$ for the monoid $\N[\frac{1}{p}].$ We fix a complete algebraic closure $C = C_K$ of $K$, a  compatible choice of $p$-power roots of unity $\zeta_{p^\infty} \hookrightarrow \mO_C$, and a compatible choice  $\varpi^{N_\infty} \hookrightarrow \mO_C$ of $p$-power roots of $\varpi$. We write 
\[\varpi^\flat = (\varpi, \varpi^{1/p}, \varpi^{1/p^2}, ...), \epsilon = (1, \zeta_p, \zeta_{p^2}, ...)  \in \mO_C^\flat\] and define elements 
\[ q = [\epsilon], \quad \mu = q - 1, \quad \xi = \mu/\varphi^{-1}(\mu), \quad \sq \xi = \varphi(\xi)\] in $\Ainf = W(\mO_C^\flat)$ as usual. We identify $(\Prism_{\mO_C}, I) = (\Ainf, \sq \xi)$ as in \cite{BS_crystal} and write $\Binf = \Ainf[1/p]$. 

\subsection{The absolute log prismatic site of $(\mO_K, \varpi^\N)$}

Let us first record some examples that will be useful later on.

\beg \label{example:BKF}   Using the quasisyntomic cover $(\mO_K, \varpi^\N) \ra (\mO_C, \varpi^{N_\infty})$, we obtain an object 
\begin{equation} \label{eq:the_Ainf_log_prism} (\Prism_{(\mO_C,N_\infty)}, I, [\varpi^{\flat}]^{N_\infty}) \cong (\Prism_{\mO_C}, I, N_\infty) 
\end{equation} in the absolute log prismatic site  $(\mO_K, \varpi^\N)^{\log}_{\Prism}$ that covers the final object of $\tu{Shv}((\mO_K, \varpi^\N)^{\log}_{\Prism})$. Evaluating a log prismatic $F$-crystal $\mE \in   \tu{Vect}^{\varphi} ((\mO_K, \varpi^\N)^{\log}_{\Prism}, \mO_{\Prism})$ on this pre-log prism corresponds to the natural pullback functor 
\[
 \tu{Vect}^{\varphi} ((\mO_K, \varpi^\N)^{\log}_{\Prism}, \mO_{\Prism}) \lra  \tu{Vect}^{\varphi} ((\mO_C, \varpi^{N_\infty})^{\log}_{\Prism}, \mO_{\Prism}).
\]
There is an intermediate cover of the final object of $\tu{Shv}((\mO_K, \varpi^\N)^{\log}_{\Prism})$ given by $(\Prism_{\mO_C}, I, [\varpi^\flat]^\N)$. Also note that the forgetful functor $(\mO_C, \varpi^{N_\infty})^{\log}_{\Prism} \ra (\mO_C)_{\Prism}$ sending a pre-log prism $(A, I, M_A)$ to its underlying prism $(A, I)$ induces an equivalence of categories 
\begin{equation} \label{eq:BKF}
     \tu{Vect}^{\varphi} ((\mO_C)_{\Prism}, \mO_{\Prism}) \isom   \tu{Vect}^{\varphi} ((\mO_C, \varpi^{N_\infty})^{\log}_{\Prism}, \mO_{\Prism}) \isom  \tu{Mod}^{\Bdrp} (\Z_p).  
\end{equation}
Here $\tu{Mod}^{\Bdrp} (\Z_p)$ denotes the category of pairs $(T, F)$ consisting of a finite free $\Z_p$-modules $T$ and a $\Bdrp$-lattice $F \subset T \otimes_{\Z_p} \Bdr$, and the second equivalence uses the classification of Breuil--Kisin--Fargues modules (see \cite[Theorem 5.2]{BS_crystal}). 
\eeg 

\beg \label{example:Acrys}
Extending scalars from (\ref{eq:the_Ainf_log_prism}) we have pre-log prisms 
\[
(\Acrys, (p), [\varpi^\flat]^{\N}) 
 \quad \tu{and} \quad (\Acrys, (p), [\varpi^\flat]^{N_\infty}) %\cong (\Prism_{(\mO_C/p, N_\infty)}, I, N_\infty)
\]
in $(\mO_K, \varpi^\N)^{\log}_{\Prism}$ that are of crystalline nature. In other words, they live in $(\mO_K/p, \varpi^\N)^{\log}_{\Prism}.$ 
\eeg 
 
\beg[The Breuil--Kisin pre-log prism] \label{example:BK} 
Let $E \in \mathfrak S := W[\![u]\!]$ denote the Eisenstein polynomial of $\varpi$ and consider the pre-log prism $(\mathfrak S, E, u^\N)$.  We will refer to this as the \textit{Breuil--Kisin pre-log prism}. There are natural maps $(\mathfrak S, E, u^{\N}) \ra (\Ainf, (\xi), [\varpi^{\flat}]^{\N}) \xrightarrow{\varphi} (\Ainf, I, [\varpi^{\flat}]^{\N})$ of pre-log prisms where the second Frobenius map is $\varphi_{\Ainf}$ on $\Ainf$ and multiplication by $p$ on the monoid. 
\eeg

\beg[The log prism $\mS$ (\`a la Kato and Breuil)] \label{example:Breuil} Let $\mS$ be the $p$-completed PD envelop with respect to the kernel of the map $W[u] \ra \mO_K$ sending $u \mapsto \varpi$.\footnote{$\mS$ is sometimes denoted by $S$ in the literature, which is slightly different from \cite{Breuil}, which denoted $\mS$ by $S_{\min}^0$. Our $\mS[1/p]$ is denoted by $S_{\min}$ in \textit{loc.cit}.} More explicitly 
\[\mS =  W[u, \frac{E(u)^n}{n!}]^{\wedge}_p.\] 
We equip $\mS$ with the pre-log structure $\N \ra \mS$ sending $1 \mapsto u$. By a slight abuse of notation, we also use $\mS$ to denote the pre-log prism $(\mS, (p), u^\N)$, which is an object in $(\mO_K, \varpi^\N)^{\log}_{\Prism}$ via the Frobenius   
\[\mO_K \cong W[u]/E(u) \xrightarrow{\varphi} \mS/p.\] By construction, $\mS$ is also an object in the ($p$-completed) log crystalline site of $(\mO_K/p, \varpi^\N)$ over $(W,(p))$. In fact, it can be regarded as an object in the strict $\delta_{\log}$-crystalline site\footnote{This is slight variant of \cite[Definition 6.4]{Koshikawa}.}  $(\mO_K/p, \varpi^\N)_{/W}^{\delta_{\tu{Crys}},\tu{str}}$ which consists of a bounded pre-log prism $(B, (p), M_B)$ over $(W, (p))$, equipped with a $p$-completed PD-ideal $J \subset B$ such that $B/J$ is classically $p$-complete, and equipped with a map 
\[(\mO_K/p, \varpi^\N) \ra (B/J, M_B)
\]
of pre-log rings that induces a strict map on the corresponding log formal schemes.   Now let us form the self-product of $\mS$ in  $(\mO_K/p, \varpi^\N)_{/W}^{\delta_{\tu{Crys}},\tu{str}}$ (or in the log crystalline site), then we get a pre-log ring 
\[
\mS^{(1)} := W[u, v, \frac{E^n}{n!}, \frac{(v/u - 1)^n}{n!}]^{\wedge}_p \cong \mS \gr{t}^{\tu{\tiny PD}}, \quad N^{(1)} \ra \mS^{(1)}
\]
where $t$ denotes the element $v/u - 1$, and 
\begin{equation} \label{eq:self_product_exactified}
N^{(1)} = \{(x, y) \in \Z^{\oplus 2} | x + y \ge 0\}
\end{equation} is the exactification of the addition map $\N^{\oplus 2} \ra \N$, and the pre-log structure is given by 
\[(x,y)\in N^{(1)} \longmapsto u^x v^y = u^{x+y} (t+1)^y \in \mS^{(1)}.\] Both arrows $p_0, p_1: \mS \rightrightarrows 
\mS^{(1)}$ are covers in  $(\mO_K/p, \varpi^\N)_{/W}^{\delta_{\tu{Crys}},\tu{str}}$ with respect to the flat topology. In particular, for any $F$-crystal $\mE \in \tu{Vect}^{\varphi} ((\mO_K/p, \varpi^\N)_{/W}^{\delta_{\tu{Crys}},\tu{str}})$ (for example, ones coming from pullbacks from $(\mO_K, \varpi^\N)^{\log}_{\Prism}$), we have an isomorphism $\beta: p_1^* \mE \isom p_0^* \mE $
that satisfies a cocycle condition when further pulled back to the triple self-product $\mS^{(2)}$ of $\mS$. Using the grading 
\begin{equation} \label{eq:grading_on_S^1}
\mS^{(1)} \cong \widehat \bigoplus_{i \ge 0} t^{[i]} \cdot \mS
\end{equation}
given by the $i^{th}$ PD-monomial $t^{[i]}$, we can further extract a $W$-linear monodromy operator $\mN = \mN_{\mS}$ on the finite free $\mS$-module %$\mS$ is local
$\hat M := \mE (\mS)$ by considering the composition
\begin{equation} \label{eq:monodromy_over_S}
\mN_{\mS}: \hat M \ra \hat M \otimes_{\mS, p_1} \mS^{(1)} \isom  \hat M \otimes_{\mS, p_2} \mS^{(1)} \xrightarrow{\tu{id}\otimes \tu{pr}^{[1]}} \hat M \otimes_{\mS} t^{[1]} \cdot \mS \cong \hat M
\end{equation}
where $\tu{pr}^{[1]}$ denotes the projection $\mS^{(1)} \ra t^{[1]} \cdot\mS$ coming from (\ref{eq:grading_on_S^1}). In particular, it induces a ``monodromy operator'' $N_{\mS}: \mS \ra \mS$ on $\mS$. It is easily checked that $N_{\mS}$ sends $u^i \mapsto i u^i$. Furthermore, the operator $\mN = \mN_{\mS}$ satisfies the following conditions: 
\begin{align} 
\label{eq:monodromy_condition1}  &  \mN (a m) = N(a)m + a \mN(m)  \quad \tu{for all } a \in \mS, m \in \hat M, \quad \tu{ and } \\
\label{eq:monodromy_condition2}  & \mN \varphi= p \varphi \mN.
\end{align}
In other words, $\hat M[1/p]$ lives in the category $\tu{Mod}^{\phi, \mN}_{/\mS[\frac{1}{p}]}$, which consists of finite free $\mS[1/p]$-modules equipped with a Frobenius $\varphi$ that is $\mS$-semilinear and invertible upon choosing an(y) basis over $\mS$, and a monodromy operator $\mN$ satisfying the analogous conditions as in (\ref{eq:monodromy_condition1}) and (\ref{eq:monodromy_condition2}). 
\eeg 

\br[The Hyodo--Kato pre-log prism] \label{remark:Breuil_HK}
There is another, more basic pre-log prism $(W, (p), 0^\N)$ which lives in $(\mO_K/p, \varpi^\N)_{/W}^{\delta_{\tu{Crys}},\tu{str}}$ (in fact, in $(k, 0^\N)_{/W}^{\delta_{\tu{Crys}},\tu{str}}$), which we will denote by $W_{\tu{HK}}$ and refer to as the \textit{Hyodo--Kato pre-log prism}. Similar to Example \ref{example:Breuil}, from an $F$-crystal on the log-crystalline site of $(k, 0^{\N})$ (or an $F$-crystal on $(k, 0^\N)_{/W}^{\delta_{\tu{Crys}},\tu{str}}$), one obtains a finite free $W$-module $M$ equipped with a monodromy operator $N$ that is $W$-linear and satisfies the condition $N \varphi = p \varphi N$ as in (\ref{eq:monodromy_condition2}). Slightly more precisely, taking the self-product of $W_{\tu{HK}}$, we obtain 
\[W_{\tu{HK}}^{(1)}= \big(W[u,v, \frac{(u/v-1)^n}{n!}]/u \big)^{\wedge}_p \cong W\gr{t}^{\tu{\tiny PD}} \] 
where $t = u/v - 1$. The monodromy operator is obtained via the same procedure as in (\ref{eq:monodromy_over_S}) (this is essentially the construction in \cite{Hyodo_Kato}). In particular, after inverting $p$ one obtains an object in $\tu{Mod}^{\phi, N}_{/K_0}$. Moreover, the map $\mS \ra W$ sending $u \mapsto 0$ and $E^{[i]} \mapsto 0$ induces a functor 
\begin{equation} \label{eq:Breuils_equivalence}\tu{Mod}^{\phi, \mN}_{/\mS[\frac{1}{p}]} \xrightarrow{ \: ``u = 0 " \: } \tu{Mod}^{\phi, N}_{/K_0}.
\end{equation}
Using iterates of the Frobenius, Breuil shows that this is in fact an equivalence of categories (slightly surprisingly, as there is no map from the Hyodo--Kato pre-log prism to the pre-log prism $\mS$). In fact, Breuil upgrades the equivalence in (\ref{eq:Breuils_equivalence}) to an equivalence between the category  $\tu{MF}_{/\mS[\frac{1}{p}]}^{\phi, \mN}$ of  \textit{filtered} $(\phi, \mN)$-modules over $\mS[1/p]$ and the category $\tu{MF}^{\phi, N}_{/K}$ of filtered $(\phi, N)$-modules (with the filtration is defined over $K$), see \cite[Section 6]{Breuil}. 
\er 

\beg[The pre-log prism $\widehat A_{\tu{st}}$]\label{example:Ast}  Kato and Breuil have introduced and studied another pre-log prism $\widehat A_{\tu{st}}$, which we recall below. Recall the objects $(\Acrys, (p), [\varpi^\flat]^\N)$ and $(\mS, (p), u^\N)$ in the absolute log prismatic site of $(\mO_K, \varpi^\N)$ from Example \ref{example:Acrys} and \ref{example:Breuil}. Both can be viewed as objects in the strict $\delta_{\log}$-crystalline site $(\mO_K/p, \varpi^\N)_{/W}^{\delta_{\tu{Crys}},\tu{str}}$. Taking their product in  $(\mO_K/p, \varpi^\N)_{/W}^{\delta_{\tu{Crys}},\tu{str}}$  (or the log crystalline site), we obtain the pre-log prism $(\widehat A_{\tu{st}}, (p), N^{(1)})$ where $N^{(1)}$ is as in (\ref{eq:self_product_exactified}) and 
\[
\widehat A_{\tu{st}} = \Acrys [u, \frac{(u/[\varpi^\flat] - 1)^n}{n!}]^{\wedge}_p = \Acrys \gr{X}^{\tu{\tiny PD}}. 
\]
Here we write $X = \frac{u}{[\varpi^\flat]}-1$. We note that the pre-log prism $(\widehat A_{\tu{st}}, (p), u^{\N})$ or $(\widehat A_{\tu{st}}, (p), [\varpi^{\flat}]^{\N})$ gives the same log prism as $(\widehat A_{\tu{st}}, (p), N^{(1)})$. By construction, $\widehat A_{\tu{st}}$ is equipped with a Frobenius operator $\varphi$ that is compatible with the Frobenius operators on $\mS$ and $\Acrys$, in particular, we have 
\[
\varphi(X) = (X+1)^p - 1. 
\]
There is also a natural $\Gal_K$-action on $\widehat A_{\tu{st}}$ that extends the $\Gal_K$-action on $\Acrys$,\footnote{Note that $\Gal_K$ acts trivially on $\mS$ (as a subring of $\widehat A_{\tu{st}}$). Also note that $\Gal_K$ does not act on the pre-log prism $(\Acrys, (p), [\varpi^{\flat}]^{\N})$ but rather on the associated log prism $(\Acrys, (p), [\varpi^{\flat}]^{\N})^a$, since the pre-log structure does not detect the roots of unity in $\mO_C$. Similarly, we have an induced $\Gal_K$-action on the associated log prism $(\widehat A_{\tu{st}}, (p), N^{(1)})$.} and we have 
\[
\gamma \cdot X = \frac{u}{\gamma \cdot [\varpi^\flat]} - 1 = \frac{[\varpi^\flat]}{\gamma \cdot [\varpi^\flat]} (X+1 ) - 1. \footnote{Our convention of the Galois action agrees with \cite{Bhatt_dR} but is slightly different from the one in \cite{Breuil}. If instead one uses $X' = [\varpi^\flat]/u - 1$ as the PD variable, then the formula would agree with \textit{loc.cit}.} 
\]
Finally,  by considering the self-product of $(\widehat A_{\tu{st}}, (p), u^{\N})$  in  $(\mO_K/p, \varpi^\N)_{/W}^{\delta_{\tu{Crys}},\tu{str}}$ as in Example \ref{example:Breuil}\footnote{Here we take the pre-log structure $u^{\N}$ on $\widehat A_{\tu{st}}$, in order to obtain the correct monodromy operator (which should be thought of as a certain differential operator along the ``$u$-direction'' and is trivial along the ``$[\varpi^\flat]$-direction'').}, we obtain a ring $
\widehat A_{\tu{st}}^{(1), {\delta_{\log}\tu{-crys}}}$ with two maps $ p_0, p_1: \widehat A_{\tu{st}} \rightrightarrows
\widehat A_{\tu{st}}^{(1), {\delta_{\log}\tu{-crys}}}$ (the cumbersome notation is intended to distinguish it from the self-product in the log prismatic site). More precisely, 
\[
\widehat A_{\tu{st}}^{(1), {\delta_{\log}\tu{-crys}}} = \Acrys [u, v, \frac{(u/[\varpi^\flat] - 1)^n}{n!}, \frac{(v/u - 1)^n}{n!}]^{\wedge}_p = \Acrys \gr{X, Y}^{\tu{\tiny PD}}
\]
where $X = u/[\varpi^\flat] - 1$ as before and $Y = v/u - 1$. The two maps $p_0, p_1$ are given by 
\[p_0: X \mapsto X, \quad p_1: X \mapsto v/[\varpi^\flat] - 1 = XY + X + Y. 
\]
Therefore, by considering the grading on $\widehat A_{\tu{st}}^{(1), {\delta_{\log}\tu{-crys}}}$ given by the PD-monomials $Y^{[i]}$ as in \ref{eq:grading_on_S^1}, we have an $\Acrys$-linear operator $N_{\tu{st}}$ on $\widehat A_{\tu{st}}$ which extends the monodromy operator on $\mS$ and sends 
\[
N_{st}(X) = X + 1. 
\]
Moreover, for an $F$-crystal $\mE \in \tu{Vect}^{\varphi} ((\mO_K/p, \varpi^\N)_{/W}^{\delta_{\tu{Crys}},\tu{str}})$, its evaluation $\sM = \mE (\widehat A_{\tu{st}})$ on $\widehat A_{\tu{st}}$ comes equipped with a monodromy operator $\mN = \mN_{\tu{st}}$ satisfying analogous conditions as in (\ref{eq:monodromy_condition1}) and (\ref{eq:monodromy_condition2}), which is $\Acrys$-linear and is compatible with the monodromy operator on $\hat M = \mE (\mS)$ coming from (\ref{eq:monodromy_over_S}) by functoriality. In other words,  along the isomorphism 
\[
\sM = \mE (\widehat A_{\tu{st}}) \cong \mE (\mS) \otimes_{\mS} \widehat A_{\tu{st}}, 
\]
induced from $(\mS, (p), u^{\N}) \ra (\widehat A_{\tu{st}}, (p), u^{\N})$, we have $\mN_{\tu{st}} = \tu{id} \otimes  N_{\tu{st}} + \mN_{\mS} \otimes \tu{id}$. On the other hand, along the isomorphism 
\[
\sM = \mE (\widehat A_{\tu{st}}) %= \mE ((\widehat A_{\tu{st}}, (p), N^{(1)})^a) 
\cong \mE ((\Acrys, (p), [\varpi^\flat]^\N)) \otimes_{\Acrys} \widehat A_{\tu{st}},
\]
induced from $(\Acrys, (p), [\varpi^\flat]^{\N}) \ra (\widehat A_{\tu{st}}, (p), N^{(1)})^a$, we have $\mN_{\tu{st}} = \tu{id} \otimes N_{\tu{st}}$. %, since the two maps from $\Acrys \ra \widehat A_{\tu{st}}^{(1)}$ are identical. 
\eeg

\br[Kato] In \cite[Section 3]{Kato_st_rep}, Kato gives a more geometric interpretation of $\widehat A_{\tu{st}}$ as  %Consider the maps of pre-log rings %$(W(k)[u], u^\N) \ra (\mO_K, \varpi^\N)$, this induces a map 
%$
%  (W(k)[u], u^\N) \ra (\mS, u^{\N}) %\ra (\Ainf, [\varpi^\flat]^{N_\infty}) 
%  \ra  (\mO_C, \varpi^{N_\infty})$, then we have  
\[
\widehat A_{\tu{st}}  \cong  R \Gamma_{\logcrys}((\mO_C, \varpi^{N_\infty})/(\mS, u^{\N})) %\cong \widehat \L \Omega_{ (\mO_C, N_\infty)/\ul{S}}. 
\]
%(see also \cite[Section 2]{Breuil}) 
Also see \cite[Remark 9.25]{Bhatt_dR} for a derived (log) de Rham cohomology interpretation of $\widehat A_{\tu{st}}$. While we  will not need these geometric interpretations, they suggest that $\widehat A_{\tu{st}}$ (and thus $\widehat B_{\tu{st}} = \widehat A_{\tu{st}} [1/p\mu]$) seems to be a more natural object to consider compared to the standard semistable period ring $B_{\tu{st}}$, which, on the other hand, is slightly simpler to use for the definition of semistable representations (see \cite[Definition 3.2, Proposition 3.3]{Breuil}).
\er 

\br \label{remark:diagram} The various pre-log prisms fit into the following diagram 
\[
\begin{tikzcd}%[row sep = 1.2em]
 & (\Prism_{\mO_C}, I, [\varpi^\flat]^{\N}) \arrow[d] \arrow[r]  & (\Prism_{\mO_C}, I,  [\varpi^\flat]^{N_\infty}) \arrow[d]  \\
(\mathfrak S, E, u^{\N}) \arrow[ru, dashed, "\iota_{\varphi}"] \arrow[d, "\varphi"]  & (\Acrys, (p), [\varpi^\flat]^{\N})  \arrow[d] \arrow[r]
& (\Acrys, (p), [\varpi^\flat]^{N_\infty}) \arrow[d] \\
(\mS, (p), u^\N)  \arrow[ru, dashed, "\iota"] \arrow[r, "\iota_{\tu{st}}"] \arrow[d]  & (\widehat A_{\tu{st}}, (p), N^{(1)}) \arrow[r] & (\widehat A_{\tu{st}}, (p), N^{(1)}_\infty)  \\
(W, (p), 0^{\N})
\end{tikzcd}
\]
where $N^{(1)}_\infty = N^{(1)} \oplus_{\N} N_\infty)$ (in other words, the bottom right square is a pushout square). The two dashed arrows are given by $\iota_{\varphi}: u \mapsto [\varpi^\flat]^p$ (resp. $\iota: u \mapsto [\varpi^\flat]$), which are both Frobenius equivariant but \textit{not} $\Gal_K$-equivariant (though they are $\Gal_{K_\infty}$-equivariant, where $K_\infty = K(\varpi^{1/p^\infty})$). On the other hand, the inclusion $\iota_{\tu{st}}$ is $\Gal_K$-equivariant %(on the associated log prisms) 
where the $\Gal_K$-action on $\mS$ is trivial, as explained in Example \ref{example:Ast}. The three squares in the diagram are commutative (but not the triangle). 
\er

\subsection{Admissibility of the \'etale realization}
Next we show that the \'etale realization functor on $F$-crystals on $(\mO_K, \varpi^{\N})^{\log}_{\Prism}$ takes values in $\Z_p$-lattices inside semistable $\Gal_K$-representations. 

\bp \label{prop:Bst_admissible}
Let $g: \tu{Vect}^{\varphi} \big((\mO_K, \varpi^\N)^{\log}_{\Prism}, \mO_\Prism  \big) \lra \tu{Vect}^{\varphi} \big((\mO_K, \varpi^\N)^{\log}_{\Prism}, \mO_\Prism[1/I_{\Prism}]^\wedge_p \big)$ be the base change functor and let 
\[
T: \tu{Vect}^{\varphi} \big((\mO_K, \varpi^\N)^{\log}_{\Prism}, \mO_\Prism  \big)  \lra \tu{Rep}_{\Z_p} (\tu{Gal}_K)
\] be the composition of the functor $g$ and the \'etale realization functor (\ref{eq:etale_realization}). Then for any $F$-crystal $\mE \in  \tu{Vect}^{\varphi} \big((\mO_K, \varpi^\N)^{\log}_{\Prism}, \mO_{\Prism} \big)$, $T(\mE)[1/p]$ is a semistable $\tu{Gal}_K$-representation. 
\ep 

\bproof  Let $\mE$ be an $F$-crystal of vector bundles on $(\mO_K, \varpi^\N)^{\log}_{\Prism}$ of rank $n$, then $T(\mE)[1/p]$ is a $\Q_p$-valued $\Gal_K$-representation of dimension $n$. 
By the equivalence (\ref{eq:BKF}) and the classification of Breuil--Kisin--Fargues modules (see, for example, \cite[Theorem 5.2]{BS_crystal}), we have a $\Gal_K$-equivariant isomorphism 
\[
\mE ((\Prism_{\mO_C}, I, [\varpi^{\flat}]^{N_\infty}))[1/\mu] \cong T(\mE) \otimes_{\Z_p} \Prism_{\mO_C} [1/\mu]
\]
when evaluating $\mE$ on the top right corner of the diagram in Remark \ref{remark:diagram} (and inverting $\mu$). Following the right vertical arrows and bottom horizontal arrows in the diagram,\footnote{In fact, we consider the corresponding ``log prisms'' associated to the pre-log prisms in order to obtain the Galois action, as remarked in Example \ref{example:Ast}.} we have natural isomorphisms
\begin{equation} \label{eq:admissible1}
T(\mE) \otimes_{\Z_p} \widehat A_{\tu{st}}[1/\mu] \cong \mE (\widehat A_{\tu{st}}) [1/\mu] \cong \mE(\mS) \otimes_{\mS} \widehat A_{\tu{st}}[1/\mu]
\end{equation}
which is equivariant for the Galois action, Frobenius and monodromy operators, where $\Gal_K$ acts trivially on $\mE(\mS)$ and the monodromy operator acts trivially on $T(\mE)$ (see the end of Example \ref{example:Ast} for the discussion on the monodromy operators; also see Example \ref{example:Breuil}). Now we observe that, after further invertible $p$ and taking the $N$-nilpotent subspace, we have 
\begin{equation} \label{eq:admissible2}
T(\mE) \otimes_{\Z_p} \Bst = \big(T(\mE) \otimes_{\Z_p} \widehat A_{\tu{st}}[1/p\mu] \big)_{\tu{N-nilp}},
\end{equation}
since $\Bst = (\widehat B_{\tu{st}} )_{\tu{N-nilp}}$ as a subring of $ \widehat B_{\tu{st}}  := \widehat A_{\tu{st}}[1/p\mu]$. By Remark \ref{remark:Breuil_HK}, we have an isomorphism  
\begin{equation} \label{eq:admissible4}
\mE(\mS) [1/p] \cong D_0
\otimes_{K_0} \mS[1/p] \end{equation} 
of $(\varphi, \mN)$-modules over $\mS[1/p]$,  
where $D_0 := \mE(W_{\tu{HK}})[1/p] = \mE ((W, (p), 0^{\N}))[1/p]$ is an $n$-dimensional $K_0$-vector space. Therefore, combining (\ref{eq:admissible1}), (\ref{eq:admissible2}), and (\ref{eq:admissible4}), we have an inclusion 
\begin{equation} \label{eq:admissible3}
D_0 \subset \big(\mE(S) \otimes_{\mS} \widehat B_{\tu{st}}\big)_{\tu{N-nilp}}^{ \Gal_K} \cong \big(T(\mE) \otimes_{\Z_p} \widehat B_{\tu{st}}\big)_{\tu{N-nilp}}^{ \Gal_K} \cong (T(\mE) \otimes_{\Z_p} \Bst )^{\Gal_K}.
\end{equation}
Therefore, $T(\mE)[1/p]$ is semistable and the inclusion in (\ref{eq:admissible3}) is in fact an equality. 
\eproof

\begin{comment}
\br[Gee--Liu] \label{remark:Gee}
Proposition \ref{prop:Bst_admissible} may also be proved as follows.  We observe that, for \textit{any} choice of compatible $p$-power roots of $\varpi$, we have a map $\iota_{[\varpi^\flat]}: \mathfrak S \ra \Ainf$ sending $u \mapsto [\varpi^\flat]^p$ (see Remark \ref{remark:diagram}), and this map is $G_{\infty} = \Gal_{K(\varpi^{1/p^\infty})}$ equivariant. Thus, we have an isomorphism 
\[
\fM \otimes_{\mathfrak S, \iota_{[\varpi^\flat]}} W(C^\flat) \cong \mE(\Prism_{\mO_C}) \otimes W(C^\flat) \cong T \otimes W(C^\flat) 
\]
which is $G_\infty$- and $\varphi$-equivariant (note that the first map is \textit{not} $\Gal_K$-equivariant),
where $T = T(\mE) $ %= (\mE(\Prism_{\mO_C}) \otimes W(C^\flat))^{\varphi=1}
denotes the \'etale realization of $\mE$. One can then show that 
\er 
\end{comment}

\br \label{remark:monodromy_on_D_0}  The proof above in fact shows that the isocrystal $D_0  = \mE(W_{\tu{HK}}) [1/p]$ agrees with $D_{\tu{st}} (T(\mE)[1/p])$ constructed by Fontaine. Moreover, the monodromy operator $N$ on $D_{\tu{st}}$ induced from $\tu{id}\otimes N$ on $T(\mE) \otimes_{\Z} \Bst$ via the inclusion (\ref{eq:admissible3}) agrees with the monodromy operator coming from the evaluation $D_0 = \mE (W_{\tu{HK}})[1/p]$ of the $F$-crystal on the Hyodo--Kato prism $W_{\tu{HK}} = (W, (p), 0^{\N})$, as constructed in Example \ref{remark:Breuil_HK}. 
\er 

\br \label{remark:cube} Gathering what we have discussed so far, we have a commutative diagram
\[
{\scalefont{0.8}
\begin{tikzcd}[column sep = -1.2cm, row sep = 0.4cm]
\tu{Vect}^{\varphi} \big((\mO_K, \varpi^\N)^{\log}_{\Prism}, \mO_\Prism  \big) \arrow[rr, "g"] \arrow[ddd, swap, "T"] \arrow[rd] &  &  \tu{Vect}^{\varphi} \big((\mO_K, \varpi^\N)^{\log}_{\Prism}, \mO_\Prism[1/I_{\Prism}]^\wedge_p \big) 
\arrow[rd] \arrow[ddd, equal, swap, near end, "T"]
%\arrow{dd}{\rotatebox{270}{$\sim$}}[swap]{T} 
\\ 
& 
\tu{Vect}^{\varphi} \big((\mO_C, N_\infty)^{\log}_{\Prism}, \mO_\Prism  \big) \arrow[rr, crossing over]   &  &  \tu{Vect}^{\varphi} \big((\mO_K, N_\infty)^{\log}_{\Prism}, \mO_\Prism[1/I_{\Prism}]^\wedge_p \big) \arrow[ddd, equal]
\\ \empty & 
\\ 
\tu{Rep}^{\tu{st}}_{\Z_p} (\Gal_K) \arrow[rd, "h"] \arrow[rr] \arrow[ruu, dashed] & &  \tu{Rep}_{\Z_p} (\Gal_K) \arrow[rd]
\\ 
& \tu{Mod}^{\Bdrp} (\Z_p) \arrow[rr] \arrow[uuu, equal, crossing over, near end, "T_\infty"] & & \tu{Mod} (\Z_p). 
\end{tikzcd}  }
\]
where $\tu{Mod} (\Z_p)$ denotes the category of finite free $\Z_p$-modules and  $\tu{Mod}^{\Bdrp} (\Z_p)$ is defined in Example \ref{example:BKF}. 
In the diagram, equal arrows are used to denote equivalence of categories (rather than equalities of any sort) to avoid messy presentation. To see that the left most square indeed commutes, we let $T = T(\mE)[1/p]$, then we have
\[h(T) = (T, D_{\tu{st}}(T[1/p])\otimes_{K^0} \Bdrp) \in \tu{Mod}^{\Bdrp} (\Z_p).\] On the other hand, 
\[T_{\infty} (\mE(\Prism_{(\mO_C, N_\infty)})) = (T, \mE(\Prism_{(\mO_C, N_\infty)})\otimes_{\Prism_{\mO_C}} \Bdrp
).\] 
To check the desired commutativity we have to identify the two $\Bdrp$-lattices in $T \otimes_{\Z_p} \Bdr$. For this, we observe that, by the proof of Proposition \ref{prop:Bst_admissible}, we have 
\[D_{\tu{st}} (T[1/p]) \otimes_{\Z_p} \Bdrp \cong \mE((W, (p), 0^{\N})) \otimes_{W} \Bdrp \cong \mE(\mS) \otimes_{\mS} \Bdrp,  
\]
where the last term identifies with $\mE(\Prism_{(\mO_C, N_\infty)})\otimes_{\Prism_{\mO_C}} \Bdrp$ by base change (see the bottom dashed arrow in the diagram in Remark \ref{remark:diagram}).  
\er

\subsection{The \'etale realization is fully faithful}
\bp 
The natural base change functor 
\begin{equation} \label{eq:base_change_fully_faithful_main}
g: \tu{Vect}^{\varphi} \big((\mO_K, \varpi^\N)^{\log}_{\Prism}, \mO_\Prism  \big) \lra \tu{Vect}^{\varphi} \big((\mO_K, \varpi^\N)^{\log}_{\Prism}, \mO_\Prism[1/I_{\Prism}]^\wedge_p \big) 
\end{equation}
is fully faithful. In particular, the \'etale realization functor  
\begin{equation} \label{eq:etale_into_st_reps}
T_{\tu{st}}: \tu{Vect}^{\varphi} \big((\mO_K, \varpi^\N)^{\log}_{\Prism}, \mO_\Prism  \big) \lra \tu{Rep}^{\tu{st}}_{\Z_p} (\tu{Gal}_K).
%\tu{Vect}^{\varphi} \big((\mO_K, \varpi^\N), \Prism_{\bullet} \big) \lra \tu{Rep}^{\tu{st}}_{\Z_p} (\tu{Gal}_K). 
\end{equation}
 from $F$-crystals on $(\mO_K, \varpi^\N)^{\log}_{\Prism}$ to $\Z_p$-lattices in semistable $\Gal_K$-representations is  fully faithful.
\ep 

\bproof 
It suffices to prove the corresponding statement for $(\mO_K, \varpi^\N)_{\qrsp}$ and $\Prism_{\bullet}$ (and $\Prism_{\bullet} [1/I]^{\wedge}_p$). For this, we take the quasisyntmic cover $(\mO_K, \N) \ra (\mO_C, N_\infty)$ and consider its ($p$-completed) Cech nerve by $(\mO_C^{(\bullet)}, N_\infty^{(\bullet)})$. We claim that the functor 
\begin{equation} \label{eq:faithfulness_on_Cech}
\tu{Vect}^{\varphi} \big((\mO_C^{(i)}, N_\infty^{(i)}), \Prism_{\bullet} \big) \lra \tu{Vect}^{\varphi} \big((\mO_C^{(i)}, N_\infty^{(i)}), \Prism_{\bullet}[1/I]^\wedge_p \big) 
\end{equation}
is faithful for each term in the Cech nerve. This follows from the logarithmic analogue of \cite[Lemma 4.11]{BS_crystal}, which in turn follows from the Hodge-Tate comparison of derived log prismatic cohomology \cite[Proposition 4.5]{logprism}. In particular, this implies that the functor (\ref{eq:base_change_fully_faithful_main}) is faithful. For fullness, we use a similar argument as in \cite{logprism} which relies on Fargues's equivalence (see Example \ref{example:BKF} and the arrow $T_\infty$ in Remark \ref{remark:cube}). Let $\mE_1, \mE_2$ be two $F$-crystals in $\tu{Vect}^{\varphi} \big((\mO_K, \varpi^\N)^{\log}_{\Prism}, \mO_\Prism  \big)$ and let $\sq \mE_1, \sq \mE_2$ be their images along 
\[
\tu{Vect}^{\varphi} \big((\mO_K, \varpi^\N)^{\log}_{\Prism}, \mO_\Prism  \big) \lra \tu{Vect}^{\varphi} \big((\mO_C, N_\infty)^{\log}_{\Prism}, \mO_\Prism  \big). 
\]
We want to show that any map $f: g(\mE_1) \ra g(\mE_2)$ comes from a map $f_0: \mE_1 \ra \mE_2$. The map $f$ gives rise to a map $T (\mE_1) \ra T (\mE_2)$ of $\Gal_K$-representations via the \'etale realization, which in turn gives a map $\sq f: \sq \mE_1 \ra \sq \mE_2$ (using the dashed arrow in Remark \ref{remark:cube}). To show that the map $\sq f$ descents to a map $f_0: \mE_1 \ra \mE_2$, it suffices to check that $\sq f$ pulls back to the same map in $
\tu{Vect}^{\varphi} \big((\mO_C^{(1)}, N_\infty^{(1)}), \Prism_{\bullet} \big)$ 
along the two projections $(\mO_C, N_\infty) \rightrightarrows (\mO_C^{(1)}, N_\infty^{(1)})$. By the faithfulness of the functors in (\ref{eq:faithfulness_on_Cech}), it suffices to check this in $\tu{Vect}^{\varphi} \big((\mO_C^{(1)}, N_\infty^{(1)}), \Prism_{\bullet}[1/I]^\wedge_p \big),$ which follows from the construction. 
\eproof

\section{Proof of the equivalence} \label{section:proof_of_equivalence} 

In this section we show the functor (\ref{eq:etale_into_st_reps}) is essentially surjective, thus completing the proof of Theorem \ref{mainthm:equivalence} of the article.   

\subsection{A construction of Berger and Kisin} \label{ss:Kisin_Berger_construction} We first describe a compatible sequence of functors 
\begin{equation} \label{eq:the_bundles_psi_n}
\psi_n: \tu{MF}_{/K}^{\phi, N} \lra  \tu{Vect}^{\varphi} \big((\mO_K, \varpi^\N)_{\qrsp}, \Prism_{\bullet} \gr{\varphi^n(I)/p}[1/p] \big) 
\end{equation}
from %the category of 
filtered $(\phi, N)$-modules to $F$-crystals of $ \Prism_{\bullet} \gr{\varphi^n(I)/p}$-vector bundles on $(\mO_K, \varpi^\N)_{\qrsp}$ (see Subsection \ref{ss:qsyn_site}), where the compatibility is with respect to the maps $\Prism_{\bullet}\gr{\varphi^n(I)/p} \ra \Prism_{\bullet}\gr{\varphi^{n-1}(I)/p}$.\footnote{which may be intuitively thought of as restricting to disks of smaller radius.} The construction is a logarithmic version of the constructions in \cite[Section 6.2]{BS_crystal} and essentially follows the recipe of Kisin in \cite[Section 1.2]{Kisin_crystal} (which is inspired by a similar construction of Berger in \cite{Berger_Robba}). 

\subsubsection{Construction of $D^{\log}$}
%Let us first describe $\psi_0$.
Let $D = (D, \varphi_D, N, \Fil^\bullet D_K) \in  \tu{MF}_{/K}^{\phi, N}$ where $D$ is a $K_0$-vector space of dimension $d$. From $D$ we can extract another filtered $\phi$-module (this time with the trivial monodromy action) as follows. Let $K_0[l_u]$ be the polynomial ring with one variable $l_u$ over $K_0$, equipped with a Frobenius map $\varphi$ extending the Frobenius on $K_0$ and sending $l_u$ to $p l_u$, and a $K_0$-linear operator $N$ acting as differential with respect to $l_u$ (one may think of $K_0 [l_u]$ formally as $K_0 [\log (u/[\varpi^\flat])] \subset \widehat A_{\tu{st}}$ from Example \ref{example:Ast}). We then define 
\[
D^{\log} := (D \otimes_{K_0} K_0 [l_u])^{N= 0}  
\]
where $N$ acts on $D \otimes_{K_0} K_0 [l_u]$ as $\tu{id}\otimes N + N \otimes \tu{id}$
Note that on $K_0 [l_u]$ we have $N \varphi = p \varphi N$, so $D^{\log}$ is a $K_0$-vector space of dimension $d$ equipped with an isomorphism $\varphi: D^{\log} \isom D^{\log}$ coming from the Frobenius $\varphi_D \otimes \varphi$ on $D \otimes_{K_0} K_0 [l_u]$. Finally, we equip $D_K \otimes_K K[l_u]$ with filtration coming from $\Fil^\bullet D_K$ and the 
$(l_u)$-adic filtration on $K[l_u],$ and equip 
\[D^{\log}_K = D^{\log} \otimes_{K_0} K \cong (D_K \otimes K[l_u])^{N = 0}\] with the subspace filtration. 
This way we obtain a functor $\tu{MF}_{/K}^{\phi, N} \ra \tu{MF}_{/K}^{\phi}$ from filtered $(\phi, N)$-modules to filtered $\phi$-modules, by sending $D \mapsto D^{\log}$. 

\subsubsection{Construction of $\psi_0$} \label{sss:psi_0}
Using $D^{\log}$, one can build the $F$-crystal $\psi_{0} (D)$ (as in \cite[Construction 6.5]{BS_crystal}): first we construct $\mM^{\tu{naive}} \in \tu{Vect}^{\varphi} \big((\mO_K, \varpi^\N)_{\qrsp}, \Prism_{\bullet} \gr{I/p}[1/p] \big) $ by 
\[
\mM^{\tu{naive}} := (\varphi^* D^{\log}) \otimes_{W} \Prism_{\bullet} \gr{I/p}, 
\]
which is equipped with a (unit) $F$-crystal structure $\varphi_{\mM^{\tu{naive}}}:\varphi^* \mM^{\tu{naive}} \isom \mM^{\tu{naive}}$ coming from the isomorphism $\varphi_D: \varphi^* D \isom D$. To obtain the desired $F$-crystal $\psi_0 (D)$ we modify $\mM^{\tu{naive}}$ along the Hodge--Tate locus $I = 0$ using the filtration on $D_K^{\log}$ defined above. In other words, using Beauville--Laszlo, we glue the vector bundles 
\begin{align*}
\mM^{\tu{naive}}[1/I] &  \in \tu{Vect}  \big((\mO_K, \varpi^\N)_{\qrsp}, \Prism_{\bullet} \gr{I/p}[1/p, 1/I] \big)   \\ 
\Fil^0 (D_K^{\log} \otimes_K \mathbb{B}_{\tu{dR}}) & \in  \tu{Vect}  \big((\mO_K, \varpi^\N)_{\qrsp}, (\Prism_{\bullet} \gr{I/p}[1/p])^{\wedge}_I  \big) 
\end{align*}
along the isomorphism 
\begin{align*} 
\Fil^0 (D_K^{\log} \otimes_K \mathbb{B}_{\tu{dR}}) [1/I] & \cong D^{\log} \otimes_{K_0} \mathbb{B}_{\tu{dR}}  \\ & \xrightarrow{\varphi^{-1}\otimes \tu{id}}  \varphi^* D^{\log} \otimes_{K_0} \mathbb{B}_{\tu{dR}} \cong \mM^{\tu{naive}} \otimes_{\Prism_\bullet \gr{I/p}} \mathbb{B}_{\tu{dR}}.
\end{align*}
In above we have identified $\mathbb{B}_{\tu{dR}}^+ \cong (\Prism_{\bullet} \gr{I/p}[1/p])^{\wedge}_I$ using (\ref{eq:identifying_bdr_sheaves}). This yields the desired vector bundle $\psi_0 (D) \in  \tu{Vect}^{\varphi} \big((\mO_K, \varpi^\N)_{\qrsp}, \Prism_{\bullet} \gr{I/p}[1/p] \big)$ equipped with an $F$- crystal structure 
\[\varphi_{\psi_0(D)}: \varphi^*(\psi_0 (D)) [1/I] \isom \psi_0 (D)[1/I]\]
coming from $\varphi_{\mM^{\tu{naive}}}$. 

\subsubsection{Construction of $\psi_n$} The $F$-crystals $\psi_n (D)$ are constructed inductively by pulling back along Frobenius and modify the poles along the Hodge--Tate divisor as in \cite[Remark 6.6]{BS_crystal}. Let us briefly recall the construction: let $\psi_{n-1} (D) \in \tu{Vect}^{\varphi} \big((\mO_K, \varpi^\N)_{\qrsp}, \Prism_{\bullet} \gr{\varphi^{n-1}(I)/p}[1/p] \big) $,   we first consider 
\begin{equation} \label{eq:extending_Kisinbundle_over_frob}
\mM^{\tu{naive}}_{n} := \psi_{n-1} (D) \otimes_{\varphi, \Prism_{\bullet}\gr{\varphi^{n-1}(I)/p}}\Prism_{\bullet}\gr{\varphi^{n}(I)/p}.
\end{equation}
Via the Frobenius $\varphi_{\psi_{n-1}(D)}$, %$\mM_n^{\tu{naive}} [1/I]$ restricts to $\psi_{n-1}(D)[1/I]$ over $\Prism_{\bullet}\gr{\varphi^{n-1}(I)/p}$. 
we have an isomorphism 
\[
\mM_n^{\tu{naive}} [1/I] \: \vline_{ \:  \Prism_{\bullet}\gr{\varphi^{n-1}(I)/p} } \isom \psi_{n-1}(D) [1/I]. 
\]
%To simplify notation, let us write $ \mathbb{B}^+_{\tu{dR}} := \big(\Prism_{\bullet}\gr{\varphi^{n}(I)/p}[1/p]\big)^{\wedge}_I [1/I] $
We then glue $\mM_n^{\tu{naive}}[1/I]$ with the $\mathbb{B}^+_{\tu{dR}}$-vector bundle $\psi_{n-1}(D)^{\wedge}_I$ along the isomorphism 
\begin{align*}
\mM_n^{\tu{naive}} \otimes_{ \Prism_{\bullet}\gr{\varphi^{n}(I)/p}} \mathbb{B}_{\tu{dR}}  
%\isom & \mM_n^{\tu{naive}} \otimes_{ \Prism_{\bullet}\gr{\varphi^{n}(I)/p}} \big(\Prism_{\bullet}\gr{\varphi^{n-1}(I)/p}[1/p]\big)^{\wedge}_I [1/I] \\
\isom  \psi_{n-1} (D) \otimes_{\Prism_{\bullet}\gr{\varphi^{n-1}(I)/p}} \mathbb{B}_{\tu{dR}}  
\mathrel{\mathop{\longleftarrow}^{ \:\sim\: }}  \psi_{n-1} (D)^{\wedge}_I \otimes_{\mathbb{B}^+_{\tu{dR}} }  \mathbb{B}_{\tu{dR}}.
\end{align*}
Here the first isomorphism uses (\ref{eq:extending_Kisinbundle_over_frob}) and the second isomorphism uses (\ref{eq:identifying_bdr_sheaves}). This finishes the construction of the compatible functors in (\ref{eq:the_bundles_psi_n}). 

%\begin{align*}
%& \mM_n^{\tu{naive}} \otimes_{ \Prism_{\bullet}\gr{\varphi^{n}(I)/p}} \big(\Prism_{\bullet}\gr{\varphi^{n}(I)/p}[1/p]\big)^{\wedge}_I [1/I] \\
%\isom & \mM_n^{\tu{naive}} \otimes_{ \Prism_{\bullet}\gr{\varphi^{n}(I)/p}} \big(\Prism_{\bullet}\gr{\varphi^{n-1}(I)/p}[1/p]\big)^{\wedge}_I [1/I] \\
%\isom & \psi_{n-1} (D)[1/I] \otimes_{\Prism_{\bullet}\gr{\varphi^{n-1}(I)/p}} \big(\Prism_{\bullet}\gr{\varphi^{n-1}(I)/p}[1/p]\big)^{\wedge}_I [1/I] \\
%\mathrel{\mathop{\longleftarrow}^{ \:\sim\: }} \:  & \psi_{n-1} (D)[1/I] \otimes_{\Prism_{\bullet}\gr{\varphi^{n}(I)/p}}  \big(\Prism_{\bullet}\gr{\varphi^{n}(I)/p}[1/p]\big)^{\wedge}_I [1/I]
%\end{align*} 

\subsubsection{Restriction to boundary}
%Let $\mathscr O \subset K_0 [\![u]\!]$ (resp. $\mathscr O_{(r,s)}\subset K_0 [\![u]\!]$, resp.  $\mathscr O_{(r,s)}^b$) denote the $K_0$-subalgebra of convergent power series on the rigid analytic open unit disc (resp. on the open annulus $\Delta_{\{r < |u| < s\}}$, resp. bounded functions in $\mathscr O_{(r, s)}$). Let $\mathscr R$ (resp. $\mathscr R^{\tu{b}}$, resp. $\mathscr R^{\tu{int}}$) denote the Robba ring (resp. the bounded Robba ring, resp. the integral Robba ring) given by $\mathscr R = \lim_{r\ra 1^{-}} \mathscr O_{(r, 1)}$ (resp. $\mathscr R^b = \lim_{r\ra 1^{-}} \mathscr O^b_{(r, 1)}$, resp. $\mathscr R^{\tu{int}} = \{\sum_{i\in \Z} a_i u^i \in \mathscr R \;|\; a_i \in W\}$). Note that $\mathscr R^{\tu{int}}$ is a DVR with $\mathscr R^{b} = \tu{Frac}(\mathscr R^{\tu{int}})$. 

Let $\sq{\mathscr R}$ (resp. $\sq{\mathscr R}^{\tu{b}}$, resp. $\sq{\mathscr R}^{\tu{int}}$) denote the (extended) Robba ring (resp. the bounded Robba ring, resp. the integral Robba ring).\footnote{If we write $B^b = \{\sum_{n \gg -\infty} [c_n] p^n \in W(C^\flat)[1/p] \tu{ s.t. } |x_n| \tu{ are bounded} \}$, equipped with norms $|\cdot|_{\rho}$ \[ \vline \sum [c_n]p^n \vline_{\rho} := p^{-\inf_{n}\{v(c_n) - n \log_p \rho\}}\]
for each $\rho \in (0, 1]$, and write $B_{(0, \rho]}$ for the completion of $B^b$ with respect to the family of norms $\{|\cdot|_\rho'\}_{0 < \rho' \le \rho}$, then $\sq{\mathscr R} := \lim_{\rho \ra 0}   B_{(0, \rho]}$. Similarly, $\sq{\mathscr R}^b = \lim_{\rho \ra 0}   B_{(0, \rho]}^b$ where $B_{(0, \rho]}^b \subset B_{(0, \rho]}$ consists of $\sum [c_n]p^n$ such that $\lim_{n \ra +\infty} |c_n| \rho^n = 0$. 
Note that a few(!) other notations are used in the literature, the standard convention seems to be 
\[\sq{\mathscr R} = \sq{B}^{\dagger}_{\tu{rig}}, \quad 
 \sq{\mathscr R}^{\tu{b}} = \sq{B}^{\dagger} =  \mE_{C^\flat}^{\dagger}, \quad \sq{\mathscr R}^{\tu{int}} = \sq A^{\dagger}=  %\mO_{\mE_{C^\flat}^{\dagger}}
\mO_{\mE_{\scaleto{C^\flat}{4pt}}^{\dagger}}. 
 \] } 
One may identify $\sq{\mathscr R}^b$ with subring of $W(C^\flat)[1/p]$\footnote{The ring $W(C^\flat)[1/p]$ is often denoted by $\mE^\dagger_{C^\flat}$ in the literature.}  \[
\sq{\mathscr R}^b = \left\{\sum_{n} [c_n] p^n \in W(C^\flat)[\frac{1}{p}]  \: \vline \: \exists \rho >0 \tu{ s.t. } |c_n|\rho^n \ra  0  \tu{ as } n \ra +\infty \right\}
\]
and 
\[\sq{\mathscr R}^{\tu{int}} = \left\{\sum_{n} [c_n] p^n \in W(C^\flat)  \: \vline \: \exists \rho >0 \tu{ s.t. } |c_n|\rho^n \ra  0  \tu{ as } n \ra +\infty \right\}.
\]
The integral Robba ring $\sq{\mathscr R}^{\tu{int}}$ is a hensalian DVR with uniformizer $p$ and $\mathscr R^{b} = \tu{Frac}(\mathscr R^{\tu{int}})$. We have natural inclusions $\Ainf \hookrightarrow \sq{\mathscr R}^{\tu{int}}$ and $\Binf = \Ainf[1/p] \hookrightarrow \sq{\mathscr R}^{b}$, and $\Ainf = \sq{\mathscr R}^{\tu{int}} \cap \Binf $ where the intersection takes place in $\sq{\mathscr R}^{b}$.  Now, the functors in (\ref{eq:the_bundles_psi_n}) give rise to compatible functors 
\begin{align*}
\sq \psi_n: \tu{MF}_{/K}^{\phi, N} \lra    \tu{Vect}^{\varphi} \big((\mO_C, \varpi^{N_\infty})_{\qrsp}, \Prism_{\bullet} \gr{\varphi^n(I)/p}[1/p] \big)    \isom   \tu{Vect}^{\varphi} (\Prism_{\mO_C} \gr{\varphi^n(I)/p}[1/p]), 
\end{align*}
%where the first arrow comes $(\mO_K, \varpi^\N) \ra (\mO_C, \varpi^{N_\infty})$ 
where the target corresponds to an $F$-crystal on $\Prism_{\mO_C}[1/p] \cong \Ainf[1/p]$, which in turn corresponds to a finite free $\Binf$-module $\mM_{\Binf}=\mM_{\infty}(\Prism_{\mO_C})$ equipped with a Frobenius isomorphism 
\[\varphi_{\mM_{\Binf}} : \varphi^* \mM_{\Binf} [1/I] \isom \mM_{\Binf} [1/I].\]
Further base changing along $\Binf \ra \sq{\mathscr R}^b$, we obtain a functor 
\begin{equation} \label{eq:functor_to_bounded_robba_ring}
\psi_{\infty, \sq{\mathscr R}^b}: \tu{MF}_{/K}^{\phi, N} \lra \tu{Mod}_{\sq{\mathscr R}^b}^{\phi}
\end{equation} from the category of filtered $(\phi, N)$-modules to $\phi$-modules over $\sq{\mathscr R}^b$ which consists of finite free $\sq{\mathscr R}^b$-modules equipped with a $\varphi_{\sq{\mathscr R}^b}$-semilinear isomorphism. From the construction above,  %the $\Binf$-module $\mM_{\Binf} = \mM_{\infty} (\Prism_{\mO_C})$ has the description 
the $\sq{\mathscr R}^b$-module $\mM_{\sq{\mathscr R}^b} = \mM_{\Binf}\otimes_{\Binf} \sq{\mathscr R}^b$ agrees with the construction of Berger in \cite{Berger_Robba} after futher base changing to $\sq{\mathscr R}$, or equivalently, with Kisin's construction $\mM(D) \otimes_{\mathscr O} \sq{\mathscr R}$ from \cite[Section 1.2]{Kisin_crystal} (up to a Frobenius twist). Therefore, by Kedlaya's theory of slopes and the observation of Berger \cite[Proposition IV.2.2]{Berger_Robba} (see also \cite[Theorem 1.3.8]{Kisin_crystal}), we know that (\ref{eq:functor_to_bounded_robba_ring}) induces a functor 
\begin{equation} \label{eq:functor_to_bounded_robba_ring_slope0}
    \psi_{\infty, \sq{\mathscr R}^b}: \tu{MF}_{/K}^{\phi, N, \tu{w.a.}} \lra \tu{Mod}_{\sq{\mathscr R}^b}^{\phi, \mu = 0}
\end{equation}
from the category of weakly admissible filtered $(\phi, N)$-modules to $\phi$-modules over $\sq{\mathscr R}^b$ of slope $0$. %Using this, one can deduce 
\bc \label{cor:chooing_BKF_module_from_Robba}
Let $\mM_{\Binf} $ be the $\phi$-module over $\Binf$ obtained from  $\{\sq \psi_{n} (D)\}_{n \ge 0}$ as above where $D$ is weakly admissible. Then there exists a Breuil--Kisin--Fargues module $(\mM_{\Ainf}, \varphi_{\mM_{\Ainf}})$ such that \[(\mM_{\Binf}, \varphi_{\mM_{\Binf}}) = (\mM_{\Ainf}, \varphi_{\mM_{\Ainf}}) \otimes_{\Ainf} \Binf.\] 
\ec

The following proof is an elaboration of the second paragraph of \cite[Section 6.4]{BS_crystal}.
\bproof 
Since $\mM_{\sq{\mathscr R}^b} = \mM_{\Binf}\otimes_{\Binf} \sq{\mathscr R}^b$ has slope $0$, there exists a $\varphi$-stable lattice $\mM_{\sq{\mathscr R}^{\tu{int}}} \subset \mM_{\sq{\mathscr R}^b}$. Then we may take $\mM_{\Ainf} := \mM_{\sq{\mathscr R}^{\tu{int}}} \cap \mM_{\Binf}$ where the intersection takes place in $\mM_{\sq{\mathscr R}^b}$. %Note that $\Ainf = \Binf \cap \sq{\mathscr R}^{\tu{int}}$
$\mM_{\Ainf}$ is a finite, torsion-free, $\varphi$-stable $\Ainf$-module. To see that it is free, we first observe there is a natural isomorphism
\begin{equation} \label{eq:inclusion_from_M_Ainf}
\mM_{\Ainf} \isom \mM_{\Ainf} \otimes_{\Ainf} \sq{\mathscr R}^{\tu{int}} \cap \mM_{\Ainf}[1/p] = \mM_{\Ainf} \otimes_{\Ainf} \sq{\mathscr R}^{\tu{int}} \cap \mM_{\Binf}.
\end{equation}
By \cite[Proposition 4.13]{BMS1} and its proof, there is a canonical inclusion $\mM_{\Ainf} \subset \mM_{\Ainf, \tu{free}} \subset \mM_{\Binf}$ from $\mM_{\Ainf}$ to a finite free $\Ainf$-module such that $ \mM_{\Ainf, \tu{free}}/\mM_{\Ainf}$ is killed by a power of the ideal $(p, I)$. Since $I$ is invertible in $\sq{\mathscr R}^{\tu{int}}$, by (\ref{eq:inclusion_from_M_Ainf}) we know that in fact $\mM_{\Ainf} = \mM_{\Ainf, \tu{free}}$. %if $(p, I)^h \cdot x \in \mM_{\Ainf},$ then $p^h \cdot x \in \mM_{\Binf}$ so $x \in \mM_{\Binf}$ and $I^h \cdot x \in \mM_{\Ainf} \otimes \sq{\mathscr R}^{\tu{int}}$ so $x \in  \mM_{\Ainf} \otimes \sq{\mathscr R}^{\tu{int}}$, so $x \in \mM_{\Ainf}$.    
\eproof 

\subsection{The Beilinson fiber sequence}
We record a logarithmic analog of \cite[Proposition 6.8]{BS_crystal} which will be used in the proof of the main result.  

\bp \label{prop:log_Beilinson_fiber} 
Let $(S, M) \in {(\mO_C, 0)}_{\qrsp}$ where $S$ is $p$-torsion free. 
We have a natural isomorphism 
\[
H^0 (\Prism_{(S, M)} \{n\}^{\varphi=1})[1/p] \isom H^0 (\Prism_{(S, M)} \gr{I/p} \{n\}^{\varphi=1})[1/p].
\]
\ep 

\bproof In the proof let us write $\ul S = (S, M)$. 
We define the Tate-twist $\Z_p(n)$ as in the non-log case:
\[
\Z_p (n) (\ul S ):= \tu{fib} \big(\Fil_{N}^n \Prism_{\ul S} \{n\} \xrightarrow{\varphi - \tu{can}} \Prism_{\ul S} \{n\}\big)
\]
where $\Fil_N^\bullet$ denotes the (derived) Nygaard filtration on log prismatic cohomology and $\varphi$ denotes the map induced from $\varphi:\Fil_{N}^n \Prism_{\ul S} \ra I^n \Prism_{\ul S}$ (see \cite[Section 5.5]{logprism}). Then set $\Q_p (n) = \Z_p (n)[1/p].$ There is a natural map\footnote{Unlike the non-log case in \cite{BS_crystal}, we do not claim that this map is an isomorphism on $H^0$. See \cite[Section 5.1]{logprism}.} 
\[
\iota: \Q_p (n) (\ul S) \lra \Prism_{\ul S}\{n\}^{\varphi = 1}[1/p] 
\]
where $ \Prism_{\ul S}\{n\}^{\varphi = 1}$ denotes $\tu{fib} (\Prism_{\ul S} \{n\} \xrightarrow{\varphi - \tu{can}} I^{-n}\otimes \Prism_{\ul S} \{n\} \isom \Prism_{\ul S} \{n\})$. By considering graded pieces of the Nygaard filtration, we have a natural isomorphism 
\begin{equation} \label{isom:rational_Tate_over_F_p}
\Q_p (n) (\ul S) \isom \Prism_{\ul S/p}\{n\}^{\varphi = 1}[1/p] \cong (\Prism_{\ul S} \{I/p\} \{n\})^{\varphi = 1} [1/p]. 
\end{equation} 
Therefore, we have natural maps 
\[
\Q_p(n) (\ul S) \lra \Prism_{\ul S}\{n\}^{\varphi = 1}[1/p] \xrightarrow{ \: f \: } (\Prism_{\ul S} \gr{I/p} \{n\})^{\varphi = 1} [1/p] \xrightarrow{\: g \: } (\Prism_{\ul S} \{I/p\} \{n\})^{\varphi = 1} [1/p]
\]
whose composition agrees with the natural map $\Q_p(n) (\ul S) \ra \Q_p (n) (\ul S/p)$ under (\ref{isom:rational_Tate_over_F_p}). The maps $f$ and $g$ above induce injections on $H^0$, thus to prove the proposition, it suffices to show that the term $H^0 ((\Prism_{\ul S} \gr{I/p} \{n\})^{\varphi = 1})[1/p]$ vanishes in $H^0(\tu{cofib}(\Q_p(n)(\ul S) \ra \Q_p(n)(\ul S/p)))$ as in the proof of \cite[Proposition 6.8]{BS_crystal}. We point out that in this argument, we only need a map $\iota$ and do not need to know that it is injective or surjective on $H^0$ (in fact, surjectivity on $H^0(\iota)$ can be deduced \textit{a posteriori}). To finish 
the proof, we need the following logarithmic version of the Beilinson fiber sequence
\[
\Q_p (n)(\ul S) \lra \Q_p (n)(\ul S/p ) \lra (\L \Omega_{\ul S/\Z_p}/\Fil_{H}^n \L \Omega_{\ul S/\Z_p} \{n\}) [1/p].
\]
This follows from left Kan extending from the log free case $(S_0, M_0) = (\Z_p[x_i, y_j]_{i \in I, j \in J}, \N^{\oplus J})$, which in turn follows from the considering the Beilinson fiber sequence for the infinite root stack %??? 
(\cite[Theorem 6.17]{AMMN}).\footnote{We identify prismatic cohomology of the infinite root stack with the log prismatic cohomology by the discussion in \cite[Section 13]{DY}. The identification of Nygaard filtrations then follow from the Hodge--Tate comparisons. See also \cite[Theorem 13.8]{DY}}
The rest of the proof follows \textit{verbatim} from the proof of \cite[Proposition 6.8]{BS_crystal}, where we identify $\L \Omega_{\ul S/\Z_p} \cong \Prism_{\ul S/p}$ using \cite[Corollary 12.18 and Corollary 12.9]{DY}. 
\eproof

Consequently, the analog of \cite[Proposition 6.10]{BS_crystal} holds in our setup. 
\bc \label{cor:extending_descent_over_boundary}
Let $\sq \mE$ be a BKF module, viewed as an $F$-crystal in $\tu{Vect}^{\varphi} ((\mO_C, N_\infty)_{\qrsp}, \Prism_\bullet)$. Suppose that $\sq \mE \gr{I/p} [1/p]$ is equipped with descent data for the cover $(\mO_K, \varpi^\N) \ra (\mO_C, N_\infty)$, in other words, suppose that there is $\mM \in \tu{Vect}^{\varphi} ((\mO_K, \varpi^{\N})_{\qrsp}, \Prism_{\bullet}\gr{I/p}[1/p])$ such that $\mM$ pulls back to $\sq \mE \gr{I/p}[1/p]$, then the descent data extends uniquely to $\sq \mE [1/p]$. 
\ec 

\bproof 
The proof of \cite[Proposition 6.10]{BS_crystal} carries over \textit{verbatim} (using Proposition \ref{prop:log_Beilinson_fiber} instead of Proposition 6.8 of \textit{loc.cit}). 
\eproof

\subsection{Descent} \label{ss:proof_of_main_equivalence} Now we can finish the proof of the main result of the article. For the setup, let $(R, P) = (\mO_C \widehat\otimes_{\mO_K} \mO_C, N_\infty\oplus_{\N} N_\infty)$ be the ($p$-completed) pushout of $(\mO_K, \N) \ra (\mO_C, N_\infty)$ along itself as in the proof of Proposition \ref{eq:base_change_fully_faithful_main}. Let $p_1, p_2: \Prism_{(\mO_C, N_\infty)} \rightrightarrows \Prism_{(R, P)}$  denote the two pullback maps. %Let $q_1, q_2$ be the image of $q \in \Prism_{\mO_C} \cong \Prism_{(\mO_C, N_\infty)}$ along $p_1, p_2$. 

\bproof[Proof of Theorem \ref{mainthm:equivalence}] With the preparations above, it remains to show that the functor $T_{\tu{st}}$ from (\ref{eq:etale_into_st_reps}) is essentially surjective, the proof of which is similar to that of \cite{BS_crystal}. Let $L \in \tu{Rep}^{\tu{st}}_{\Z_p} (\Gal_K)$. Let $D = (D, \varphi_D, N, \Fil^\bullet D_K)$ be the weakly admissible filtered $(\phi, N)$-module attached to $L[1/p]$. By Corollary \ref{cor:extending_descent_over_boundary} and its proof, there is a Breuil--Kisin--Fargues module $\sq \mE'$ corresponding to a choice of a $\varphi$-stable lattice $\mM_{\sq{\mathscr R}^{\tu{int}}}$ inside the $\varphi$-module  $\mM_{\sq{\mathscr R}^{b}} =  \psi_{\infty, \sq{\mathscr R}^b} (D)$ over $\sq{\mathscr R}^{b}$, which has slope $0$ (see (\ref{eq:functor_to_bounded_robba_ring})). From the construction, the base change 
\[\sq \mE' \otimes_{\Prism_{\mO_C}} \Prism_{\mO_C}\gr{I/p}[1/p]\] agrees with the base change of the $ \Prism_{\bullet} \gr{I/p}[1/p]$-vector bundle $\psi_0 (D)$ from $(\mO_K, \varpi^\N)_{\qrsp}$ to $(\mO_C, N_{\infty})_{\qrsp}$. Thus, by Corollary \ref{cor:extending_descent_over_boundary}, $\sq \mE'[1/p]$ naturally descents to an object 
\[\mM' \in \tu{Vect}^{\varphi} ((\mO_K, \varpi^{\N})_{\qrsp}, \Prism_{\bullet} [1/p]),\] which extends the compatible system $\psi_n (D)$ constructed in Section \ref{ss:Kisin_Berger_construction}.   %as in the proof of \cite{BS_crystal}
By considering the $\varphi$-compatible isomorphisms 
\[
\mM' (\Prism_{\mO_C}) \otimes_{\Prism_{\mO_C}} W(C^\flat) \cong \sq \mE'\otimes_{\Prism_{\mO_C}} W(C^\flat)[1/p] \cong T[1/p] \otimes_{\Q_p} W(C^\flat)[1/p]
\]
where $T = (\sq \mE' \otimes W (C^\flat))^{\varphi = 1}$, we obtain a $\Gal_K$-equivariant structure on $T[1/p]$.  After picking a $\Z_p$-lattice $T' \subset T[1/p]$ and modifying the vector bundle $\sq \mE'$ along $\{p = 0\} \subset  \spec \Prism_{\mO_C} - \{x\}$ if necessary (here $x$ denotes the closed point), we may arrange so that $T$ itself is $\Gal_K$-stable. Now, the $\Z_p$-lattice $T$ corresponds to an object  $\mE'_{\ett} \in \tu{Vect}^{\varphi} ((\mO_K, \varpi^{\N})_{\qrsp}, \Prism_{\bullet} [1/I]^{\wedge}_p)$ by Proposition \ref{prop:etale_realization}, which pulls back to $\sq \mE' \otimes_{\Prism_{\mO_C}} W(C^\flat)$ along $(\mO_C, N_\infty)_{\qrsp} \ra (\mO_K, \varpi^\N)_{\qrsp}$. In other words, the descent isomorphism 
\[
\sq \alpha:
p_1^* \sq \mE' [1/I]^{\wedge}_p [1/p] \isom p_2^* \sq \mE' [1/I]^{\wedge}_p [1/p]
\]
restricts to an isomorphism  
\[
p_1^* \sq \mE' [1/I]^{\wedge}_p   \isom p_2^* \sq \mE' [1/I]^{\wedge}_p 
\]
(which satisfies the cocycle condition). By Beaville--Laszlo glueing, this shows that the descent data on $\sq \mE' [1/p]$ extends to $\spec \Prism_{\mO_C} - \{x\}$, thus to $\spec \Prism_{\mO_C}$ (using \cite[Lemma 4.6]{BMS1}). Therefore, $\sq \mE'$ descents to an $F$-crystal $\mE' \in \tu{Vect}^{\varphi} ((\mO_K, \varpi^{\N})_{\qrsp}, \Prism_{\bullet})$ satisfying $\mE' [1/p] = \mM'$ and $T_{\tu{st}} (\mE') = T$.  Next we observe that $T[1/p] \cong L[1/p]$ since $\mE' [1/p] \otimes_{\Prism_\bullet} \Prism_{\bullet} \gr{\varphi^n(I)/p}$ agrees $\psi_n (D)$, thus we may modify the vector bundle $\sq \mE'$ using the lattice $L \subset T[1/p]$ to obtain a Breuil--Kisin--Fargues module $\sq \mE$ that corresponds to $L$. To finish the argument, we simply run the above descent argument again to obtain the desired $F$-crystal $\mE$ that recovers the lattice $L = T_{\tu{st}} (\mE)$. 
\eproof

\br 
From our proof it is evident that the equivalence $T_{\tu{st}}$ is compatible with the equivalence from \cite{BS_crystal}, via the pullbank functor of $F$-crystals from the absolute prismatic site $(\mO_{K})_{\Prism}$ to $(\mO_K, \varpi^\N)^{\log}_{\Prism}$  and the fully faithful embedding $ \tu{Rep}^{\tu{crys}}_{\Z_p} (\tu{Gal}_K) \ra  \tu{Rep}^{\tu{st}}_{\Z_p} (\tu{Gal}_K)$. 
\er 
 
%\subsection{Comparison to the strict site}

\section{Evaluation on the Breuil--Kisin prism $\mathfrak S$} \label{section:BK}

In this section, we give some descriptions of $\Z_p$-lattices in semistable $\Gal_K$ representations using (semi)linear algebra  data by evaluating the corresponding $F$-crystal on certain log prisms. We start by considering the evaluation on the Breuil--Kisin pre-log prism (Example \ref{example:BK}), which leads to a functor 
\[
\tu{Rep}^{\tu{st}}_{\Z_p} (\tu{Gal}_K) \lra  \tu{Vect}^{\varphi} (\mathfrak S),
\]
where $\tu{Vect}^{\varphi} (\mathfrak S)$ denotes $F$-crystals over $\mathfrak S$, in other words, the category of finite free $\mathfrak S$-module $\mathfrak M$ equipped with a Frobenius isomorphism  $\varphi:  \varphi_{\mathfrak S}^* \mathfrak M [1/E(u)] \isom \mathfrak M  [1/E(u)]$. However, the functor above is not fully faithful (since the category $\tu{Vect}^{\varphi} (\mathfrak S)$ embeds fully faithfully into $\tu{Rep}_{\Z_p} (\Gal_{K_\infty})$ by \cite[Proposition 2.1.12]{Kisin_crystal} or \cite[Theorem 7.2]{BS_crystal} but the composition $
\tu{Rep}^{\tu{st}}_{\Z_p} (\tu{Gal}_K) \ra 
\tu{Rep}_{\Z_p} (\tu{Gal}_{K_\infty})$ fails to be fully faithful). Therefore, we need to enlarge the target category $\tu{Vect}^{\varphi}(\mathfrak S)$ with additional data in order to obtain a useful description of $\Z_p$-lattices in semistable Galois representations. %To achieve this, one should consider either a certain monodromy operator or a certain Galois action.  

\subsection{Breuil--Kisin modules with integral monodromy operators} \label{ss:BK_mod_with_int_monodromy} One way to enlarge the category $\tu{Vect}^{\varphi} (\mathfrak S)$ is to consider a certain (integral) monodromy operator. Let  $\tu{Mod}^{\phi}_{/\mathfrak S}$ denote the full-subcategory of $\tu{Vect}^{\varphi} (\mathfrak S)$ consisting of pairs $(\fM, \varphi)$ such that $\varphi$ takes  $\varphi_{\mathfrak S}^* \fM$ into $\fM$. Equivalently, objects in $\tu{Mod}^{\phi}_{/\mathfrak S}$ are finite free $\mathfrak S$-module $\mathfrak M$, equipped with a $\varphi_{\mathfrak S}$-semilinear Frobenius $\varphi: \mathfrak M \ra \mathfrak M$ such that the cokernel of $\varphi^* \fM \ra \fM$ is killed by a power of $E(u)$. We refer to such objects as \textit{Breuil--Kisin modules} (sometimes they are called \textit{effective} Breuil--Kisin modules). 

\bd \label{def:BK_mod_int_monodromy}
Let $\tu{Vect}^{\varphi, N_{\tu{int}}} (\mathfrak S)$ (resp. $\tu{Mod}^{\phi, N_{\tu{int}}}_{/\mathfrak S}$) denote the category of objects $(\fM, \varphi)$  in $\tu{Vect}^{\varphi} (\mathfrak S)$ (resp. in $\tu{Mod}^{\phi}_{/\mathfrak S}$) equipped with a $W$-linear endomorphism $N: \fM /u \fM \ra \fM /u\fM$ satisfying $N \varphi = p \varphi N$.  
\ed 

An object in $\tu{Mod}^{\phi, N_{\tu{int}}}_{/\mathfrak S}$ is called a $(\phi, N)$-\textit{module over} $\mathfrak S$.

\br \label{remark:Kisin_category}
We would have denoted the category $\tu{Mod}^{\phi, N_{\tu{int}}}_{/\mathfrak S}$ by $\tu{Mod}^{\phi, N}_{/\mathfrak S}$, except the latter notation was already used by Kisin in \cite{Kisin_crystal} to denote the category of Breuil--Kisin modules $(\fM, \varphi)$ equipped with a \textit{rational} monodromy operator $N_{\Q}: \fM /u \fM \otimes \Q_p \ra \fM /u \fM \otimes \Q_p$ satisfying $N \varphi = p \varphi N$. We will still denote Kisin's category by  $\tu{Mod}^{\phi, N}_{/\mathfrak S}$ and refer to its object by $(\phi, N_{\Q})$-modules over $\mathfrak S$.  
\er  
\begin{construction}
Let $\mE \in \tu{Vect}^{\varphi} \big((\mO_K, \varpi^\N)^{\log}_{\Prism}, \mO_\Prism  \big)$ be an $F$-crystal on $(\mO_K, \varpi^{\N}). $
The evaluation functor $\mE(\mathfrak S) := \mE ((\mathfrak S, (E(u)), u^{\N}))$ on the Breuil--Kisin pre-log prism and the isomorphism 
\[
\mE(\mathfrak S) \otimes_{\mathfrak S} \mathfrak S/u \cong \mE (W_{\tu{HK}})
\]
of finite free $W$-modules gives rise to a functor from $ \tu{Vect}^{\varphi} \big((\mO_K, \varpi^\N)^{\log}_{\Prism}, \mO_\Prism  \big)$ to $\tu{Vect}^{\varphi, N_{\tu{int}}} (\mathfrak S)$. 
% $\tu{Mod}^{\phi, N_{\tu{int}}}_{/\mathfrak S}$. 
Composing with the inverse of the equivalence (\ref{eq:etale_into_st_reps}), we obtain a  functor 
\begin{equation} \label{eq:the_functor_D_sigma}
D_{\mathfrak S}: \tu{Rep}^{\tu{st}}_{\Z_p} (\tu{Gal}_K) \lra \tu{Vect}^{\varphi, N_{\tu{int}}} (\mathfrak S)
%\tu{Mod}^{\phi, N_{\tu{int}}}_{/\mathfrak S}
\end{equation}
from $\Z_p$-lattices in semistable Galois representations to the category of $(\phi, N)$-modules over $\mathfrak S$. 
\end{construction}

\br[Hodge--Tate weights and $E$-height] \label{remark:HT_weights} Let $\mE$ be as above. 
Unwinding the construction in Subsection \ref{ss:Kisin_Berger_construction}, we see that $\fM = \mE (\mathfrak S)$ is an (effective) Breuil--Kisin module if and only if the corresponding filtered $(\phi, N)$-module $D$ is effective, i.e., $\Fil^0 D_K = D_K$. In this case, if $\Fil^h D_K = 0$ for some $h \in \Z_{\ge 0}$, then $\fM$ has finite $E$-height equal to $h$ in the sense that $E^h$ kills $\tu{coker} (\varphi: \varphi^* \fM \ra \fM)$ (see \cite[Lemma 1.2.2]{Kisin_crystal}). 
\er 

\bt  \label{thm:D_Sigma_fully_faithful}
The functor $D_{\mathfrak S}$ in (\ref{eq:the_functor_D_sigma}) is fully faithful. In particular, %by Remark \ref{remark:HT_weights}, 
there is a  natural fully faithful functor 
\[
D_{\mathfrak S}: \tu{Rep}^{\tu{st}, \ge 0}_{\Z_p} (\tu{Gal}_K) \lra \tu{Mod}^{\phi, N_{\tu{int}}}_{/\mathfrak S}.
\]
We describe the essential image of $D_{\mathfrak S}$ in the end of Subsection \ref{ss:BK_galois} (see Remark \ref{remark:essential_image_D_sigma}). 
\et 

\bproof The strategy of the proof is similar to that of \cite[Theorem 7.9]{BS_crystal}. For completeness let us recall the argument. 
It suffices to prove the second statement (with non-negative Hodge--Tate weights). We first pass to the isogeny categories and consider the composition of the natural functors 
\begin{multline} \label{eq:D_st}
\quad \mD: \tu{Rep}^{\tu{st}, \ge 0}_{\Z_p} (\tu{Gal}_K) \otimes \Q_p \xrightarrow{D_{\mathfrak S} \otimes \Q_p} \tu{Mod}^{\phi, N_{\tu{int}}}_{/\mathfrak S} \otimes \Q_p \\
\xrightarrow{-\otimes_{\mathfrak S} \mathfrak S \gr{I/p}[1/p]}  \tu{Mod}^{\phi, N}_{/\mathfrak S\gr{I/p}[1/p]}    \xrightarrow{\:\: \mD_0 \:\: } \tu{MF}_{/K}^{\phi, N} \qquad 
\end{multline} 
where the functor $\mD_0$  constructed in a similar way as in the proof of \cite[Theorem 7.9]{BS_crystal}. More precisely, the underlying isocrystal $\mD_0 (\fM)$ is given by $\fM \otimes_{\mathfrak S} K_0$ with its induced Frobenius, the filtration on $\mD_0 (\fM)_K$ is built using the $\varphi$-action on $\fM$ as in \textit{loc.cit.} (or equivalently, by \cite[1.2.5]{Kisin_crystal}), and moreover one keeps track of the (rational) monodromy operator $N_{\Q}$ on $\mD_0 (\fM)$. Unwinding the constructions in Subsection \ref{ss:Kisin_Berger_construction} and using \cite[Proposition 1.2.8]{Kisin_crystal}, we know that $\mD_0$ has a right inverse given by $\psi_0$ from Construction \ref{sss:psi_0} (composed with the evaluating on the Breuil--Kisin pre-log prism $\mathfrak S$).  We claim that the composition $\mD$ of the functors in (\ref{eq:D_st}) agrees with Fontaine's $D_{\tu{st}}$ functor. But this can be checked in the category $\tu{Mod}^{\phi, N}_{/\mathfrak S\gr{I/p}[1/p]}$ after composing with the functor $\psi_0$, which follows from the (constructive) proof of Theorem \ref{mainthm:equivalence} in Subsection \ref{ss:proof_of_main_equivalence} (see also Remark \ref{remark:monodromy_on_D_0}). Therefore, we know that $\mD$ is fully faithful. 
By the same argument of \cite[Theorem 7.9]{BS_crystal} (using the proof of \cite[Lemma 1.2.6]{Kisin_crystal}), we know that the functor $\mD_0$ is faithful, and thus $\mD_{\mathfrak S} \otimes \Q_p$ is fully faithful. To finish the proof, let $L, L' \in \tu{Rep}^{\tu{st}, \ge 0}_{\Z_p} (\tu{Gal}_K) $ and let $f: D_{\mathfrak S} (L) \ra D_{\mathfrak S} (L')$ be a map of $(\phi, N)$-modules, then for some $n \ge 0$ we have $p^n f = D_{\mathfrak S} (g)$ for some map $g: L \ra L'$ of $\Gal_K$-representations. Now since the forgetful functor (forgetting the monodromy operator) gives rise to 
\[\tu{Mod}^{\phi, N_{\tu{int}}}_{/\mathfrak S} \ra \tu{Mod}^{\phi}_{/\mathfrak S} \ra \tu{Vect}^{\phi} ({\mathfrak S[1/E]^\wedge_p}) \cong \tu{Rep}_{\Z_p} (\Gal_{K_\infty})
\]
that agrees with the forgetful functor $\tu{Rep}^{\tu{st}, \ge 0}_{\Z_p} (\tu{Gal}_K) \ra  \tu{Rep}_{\Z_p} (\Gal_{K_\infty})$, we know that $g$ is disivible by $p^n$ by considering $D_{\mathfrak S} (g) = p^n f$ in $\Z_p$-lattices of $\Gal_{K_\infty}$-representations. 
\eproof 

\br %[Compatibility with the crystalline case]
From the construction and Remark \ref{remark:monodromy_on_D_0}, the fully faithful embedding $D_{\mathfrak S}$ above is compatible with the crystalline case in  \cite[Theorem 7.9]{BS_crystal}, via the functor $ \tu{Mod}^{\phi}_{/\mathfrak S} \ra  \tu{Mod}^{\phi, N_{\tu{int}}}_{/\mathfrak S}$ which equips a Breuil--Kisin module $(\fM, \varphi)$ with the monodromy operator $N = 0$ on $\fM/u$. 
\er 

\br \label{remark:nabla}  Let $\mathscr O \subset K_0 [\![u]\!]$  denote the $K_0$-algebra of convergent power series on the rigid analytic open unit disc and let $\lambda = \prod_{n \ge 0} \varphi (\frac{ (E(u))}{  (E(0))}) \in \mathscr O$. Let  $\fM$  be a $(\phi, N)$- or  $(\phi, N_{\Q})$-module over $\mathfrak S$, then its base change $\fM_{\mathscr O} = \fM \otimes_{\mathfrak S} \mathscr O$ gives rise to a $(\phi, N)$-module over $\mathscr O$ in $\tu{Mod}^{\phi, N}_{/\mathscr O}$  with $\varphi_{\fM_{\mathscr O}} = \varphi_{\fM} \otimes \varphi_{\mathscr O}$ (see \cite[1.2.4]{Kisin_crystal}). Moreover, by \cite[Lemma 1.3.10]{Kisin_crystal}, there exists a canonical operator $N_{\nabla}$ on $\fM_{\mathscr O} [\frac{1}{\lambda}]$ satisfying $N_{\nabla} \varphi = (p/E(0)) E(u) \varphi N_{\nabla}$ and $N_{\nabla}|_{u = 0} = N_{\Q}$. Let $\tu{Mod}^{\phi, N_{\tu{int}}, \nabla}_{/\mathfrak S}$ (resp. $\tu{Mod}^{\phi, N, \nabla}_{/\mathfrak S}$) denote the full subcategory of $\tu{Mod}^{\phi, N_{\tu{int}}}_{/\mathfrak S}$ (resp. $\tu{Mod}^{\phi, N}_{/\mathfrak S}$) spanned by  $(\phi, N)$- or  $(\phi, N_{\Q})$-modules $\fM$ such that $N_{\nabla}$ preserves $\fM_{\mathscr O} \subset \fM_{\mathscr O}[\frac{1}{\lambda}]$. By the proof of Theorem \ref{thm:D_Sigma_fully_faithful}, we know that the functor $D_{\mathfrak S}$ thereof in fact factors as fully faithful functors 
\[
D_{\mathfrak S}: \tu{Rep}^{\tu{st}, \ge 0}_{\Z_p} (\tu{Gal}_K) \lra \tu{Mod}^{\phi, N_{\tu{int}}, \nabla}_{/\mathfrak S}  \lra \tu{Mod}^{\phi, N, \nabla}_{/\mathfrak S},
\]
which becomes equivalences after passing to the isogeny categories. In particular, we have 
\[\tu{Mod}^{\phi, N_{\tu{int}}, \nabla}_{/\mathfrak S}\otimes \Q_p \isom \tu{Mod}^{\phi, N, \nabla}_{/\mathfrak S} \otimes \Q_p.\]
This partially answers a question raised in \cite{Kisin_crystal}, though we still do not know whether the category $\tu{Mod}^{\phi, N_{\tu{int}}}_{/\mathfrak S}\otimes \Q_p$ is equivalent to $\tu{Mod}^{\phi, N}_{/\mathfrak S}\otimes \Q_p$. 
\er

%\subsection{Galois representations of finite $E(u)$-height}

\subsection{Breuil--Kisin modules with Galois actions} \label{ss:BK_galois}  Next we consider Breuil--Kisin modules equipped with additional data coming from certain Galois actions. 
Recall from Section \ref{section:F_crystals} that we have fixed a compatible system of $\zeta_{p^\infty}$ and $\varpi^{1/p^\infty}$ in $\mO_C$, 
%$\varpi^\flat = (\varpi, \varpi^{1/p}, ...) \in \mO_C^\flat$, 
which in particular determines a cyclotomic character $\chi$ given by $\gamma \zeta_{p^\infty} = \zeta_{p^\infty}^{\chi (\gamma)}$ and a map $\iota = \iota_{\varpi^\flat}: \mathfrak S \ra \Ainf$ sending $u \mapsto [\varpi^\flat]$. Let $\iota_\varphi: \mathfrak S \ra \Ainf$ denote the Frobenius twisted composition $\varphi \circ \iota: \mathfrak S \ra \Ainf$.  Using this map, we view a finite free $\mathfrak S$-module $\fM$ as a submodule of $\fM \otimes_{\mathfrak S, \iota_{\varphi}} \Ainf$. As in Remark \ref{remark:diagram}, let $K_{\infty} = K (\varpi^{1/p^\infty})$. Let $L =K_\infty(\zeta_{p^\infty})$, and let $\hat G = \Gal(L/K)$. Pick a topological generator $\tau \in \Gal(L/K(\zeta_{p^\infty}))$, then $\hat G$ (resp. $\Gal_K$) is topologically generated by $\Gal (L/K_\infty)$ (resp. $\Gal_{K_\infty}$) and $\tau$.\footnote{In fact, if $p \ne 2$, then we have a natural isomorphism $\hat G \cong 
 \Gal(L/K(\zeta_{p^\infty})) \rtimes \Gal (L/K_\infty)
 \cong \Z_p (1) \rtimes \Gal(L/K_{\infty})$. See, for example, \cite[Lemma 5.1.2]{Tong_Liu_compositio}.} We may (and will) choose $\tau$ such that $\tau ([\varpi^\flat]) = [\epsilon]\cdot [\varpi^\flat].$ We now give precise definitions of the categories that appear in the diagram from the introduction.

\bd
Let $\tu{Mod}_{/\mathfrak S}^{\phi, \tau}$ denote the category of Breuil--Kisin modules $(\fM, \varphi) \in \tu{Mod}^{\varphi}_{/ \mathfrak S}$ equipped with a $\varphi$-commuting continuous $\Gal_K$-action on $\fM_{\Ainf} := \fM \otimes_{\mathfrak S, \iota_{\varphi}} \Ainf$ such that $\Gal_{K_{\infty}}$ acts trivially on $\fM$ and the action of $\tau$ satisfies the condition that for each $x \in \fM$, we have 
\begin{equation} \label{eq:condition_of_tau_action}
(\tau-1) x \in  
 \fM_{\Ainf} \otimes_{\Ainf}  ([\epsilon] - 1) \Acrys  = \fM \otimes_{\mathfrak S, \iota_{\varphi}} \mu \Acrys
\end{equation}
\ed 

\br 
Equivalently, one may consider (as in \cite{Caruso}) the category of Breuil--Kisin modules $(\fM, \varphi)$ equipped with a Frobenius-commuting automorphism $\tau: \fM_{\Ainf} \ra \fM_{\Ainf}$ satisfying the condition (\ref{eq:condition_of_tau_action}) above and the condition that 
\[
(\tu{id} \otimes \gamma) \circ \tau (x) = \tau^{\chi(\gamma)} (x)
\]
for all $x\in \fM$, and all $\gamma \in \Gal_{K_\infty}$ with $\chi(\gamma) \in \N$. %This description is closer to the categories considered in \cite{Caruso}. 
\er 

There are several other closely related categories.  For the setup, let $\mathfrak S_{\log}^{(1)}$ denote the self-product of the Breuil--Kisin pre-log prism $\mathfrak S$ in $(\mO_K, \varpi^\N)_{\Prism}^{\log}$. More explicitly, 
\[ \mathfrak S_{\log}^{(1)} \cong W[\![u, v]\!]\{\frac{v/u - 1}{E}\}^{\delta, \wedge} \cong W[\![u, \mathfrak y]\!]\{\frac{\mathfrak y}{E}\}^{\delta, \wedge}\] where ${\wedge}$ denotes the derived $(p, E)$-completion (see \cite[Section 2.3]{DuLiu}). Let $\iota_1, \iota_2$ denote the two maps $\mathfrak S \ra \mathfrak S_{\log}^{(1)}$ sending $u \mapsto u$ and $u \mapsto (\mathfrak y + 1) u$ respectively. There is a $\Gal_K$-equivariant embedding $\iota^{(1)}: \mathfrak S_{\log}^{(1)} \hookrightarrow \Ainf$ sending $u \mapsto [\varpi^\flat]$ and $\mathfrak y \mapsto \mu = [\epsilon] - 1.$ Let $\iota_{\varphi}^{(1)}$ denote the Frobenius twisted embedding $\varphi \circ \iota^{(1)}: \mathfrak S_{\log}^{(1)} \ra \Acrys$. Let $I^+$ denote the kernel of the composition $\mathfrak S_{\log}^{(1)} \ra \Ainf \ra W(\cl k)$ and let $I_{\mathfrak y}^+ \subset I^+$ denote the pre-image $(\iota_\varphi^{(1)})^{-1} (\mu \Acrys)$. 

% where the second map is induced from the projection of $\mO_C^\flat$ onto its residue field $\cl k$. %$\mO_C^\flat \ra \mO_C^\flat/\fm^\flat = \cl k$

\bd \label{definitions:other_similar_BK_modules}
Let $\tu{Mod}_{/\mathfrak S}^{\phi, \tau_{\tu{st}}}$ (resp.  $\tu{Mod}_{/\mathfrak S}^{\phi, \hat G}$, resp. $\tu{Mod}_{/\mathfrak S}^{\phi, G_K}$) denote the category of objects $(\fM, \varphi) \in \tu{Mod}^{\varphi}_{/ \mathfrak S}$ equipped with a continuous $\Gal_K$-action on $\fM \otimes_{\mathfrak S, \iota_1} \mathfrak S_{\log}^{(1)}$ (resp. on $\fM \otimes_{\mathfrak S, \iota_1} \mathfrak S_{\log}^{(1)}$, resp. on $\fM \otimes_{\mathfrak S, \iota_{\varphi}} \Ainf$) such that the action $\Gal_K$ commutes with $\varphi$ and satisfies the following conditions
\be 
\item $\Gal_{K_{\infty}}$ acts trivially on the submodule $\fM$; and  
\item the action of $\tau - 1$ sends each element $x \in \fM$ into $\fM \otimes_{\mathfrak S} I^+_{\mathfrak y} \mathfrak S_{\log}^{(1)}$ (resp. into $ \fM \otimes_{\mathfrak S} I^+ \mathfrak S_{\log}^{(1)}$; resp. into $\fM \otimes_{\mathfrak S, \iota_\varphi} W([\fm^\flat]) \Ainf$ ).
%\begin{align*}
%& (\tau-1) x \in    \fM \otimes_{\mathfrak S} \mathfrak y \mathfrak S_{\log}^{(1)}   \\  \big( \tu{resp. } \:  & (\tau-1) x \in    \fM \otimes_{\mathfrak S} I^+ \mathfrak S_{\log}^{(1)} \\ \big( \tu{resp. } \:  & (\tau-1) x \in   \fM \otimes_{\mathfrak S} W([\fm^\flat]) \Ainf  \big). \end{align*}
\ee 
\ed
From the definitions, we have fully faithful embeddings of categories 
\begin{equation} \label{diagram:several_BK_categories}
\begin{tikzcd}[row sep = 1.2em, column sep = 1.5em]
\tu{Mod}_{/\mathfrak S}^{\phi, \tau_{\tu{st}}} \arrow[d, hook] \arrow[r, hook] & \tu{Mod}_{/\mathfrak S}^{\phi, \tau} \arrow[d, hook] \\
\tu{Mod}_{/\mathfrak S}^{\phi, \hat G} \arrow[r, hook] & 
\tu{Mod}_{/\mathfrak S}^{\phi, G_K} 
\end{tikzcd}
\end{equation}
where the vertical arrows correspond to slightly relaxing the condition on the image of $(\tau - 1) \fM$, while the horizontal arrows correspond to allowing the $\Gal_K$-action on $\fM$ to be defined over the larger ring $\Ainf$ instead of $ \mathfrak S^{(1)}_{\log}. $ From Theorem \ref{mainthm:equivalence}, we can deduce the following
\bt  \label{thm:equivalence_with_phi_tau}
There are natural equivalences of categories 
\begin{equation} \label{eq:equivalence_for_phi_tau}
\tu{Rep}^{\tu{st}, \ge 0}_{\Z_p} (\tu{Gal}_K) \lra   \tu{Mod}^{\phi, \tau_{\tu{st}}}_{/\mathfrak S}  \lra   \tu{Mod}_{/\mathfrak S}^{\phi, \tau}.
\end{equation}
Moreover, the essential image of $\tu{Rep}^{\tu{crys}, \ge 0}_{\Z_p} (\tu{Gal}_K)$ are Breuil--Kisin modules $\fM$ in $ \tu{Mod}^{\phi, \tau_{\tu{st}}}_{/\mathfrak S} $ (resp. $\tu{Mod}_{/\mathfrak S}^{\phi, \tau}$) such that the action of $\tau$ satisfies the more strigent condition that for each $x \in \fM$,
\[
(\tau - 1) x \in \fM \otimes_{\mathfrak S}  I_{\mathfrak y, \tu{crys}}^+  \qquad (\tu{resp. } (\tau - 1) x \in \fM \otimes_{\mathfrak S}  I_{\tu{crys}}^+ ),
\]
where $ I_{\tu{crys}}^+ := \mu  \cdot \ker (\Acrys \ra W(\cl k)) \subset \Acrys$ and $ I_{\mathfrak y, \tu{crys}}^+ := (\iota_\varphi^{(1)})^{-1} ( I_{\tu{crys}}^+).$
\et

\bproof 
Let us first describe the functor. It suffices to construct a functor from  $
\tu{Rep}^{\tu{st}, \ge 0}_{\Z_p} (\tu{Gal}_K)$ to $\tu{Mod}^{\phi, \tau_{\tu{st}}}_{/\mathfrak S}$. Let $L \in \tu{Rep}^{\tu{st}, \ge 0}_{\Z_p} (\tu{Gal}_K)$ and let $\mE$ be the corresponding log prismatic $F$-crystal over $(\mO_K, \varpi^\N)$. As in Subsection \ref{ss:BK_mod_with_int_monodromy}, evaluation on the Breuil--Kisin pre-log prism $\mathfrak S$ gives a Breuil--Kisin module $\fM = \mE(\mathfrak S) \in \tu{Mod}^{\phi}_{/\mathfrak S}$. Using the isomorphism $\alpha_1:  \mE(\mathfrak S_{\log}^{(1)}) \isom \fM \otimes_{\mathfrak S, \iota_1} \mathfrak S_{\log}^{(1)}$, we can equip the target with the desired $\Gal_K$-action coming from the Galois action on $ \mE(\mathfrak S_{\log}^{(1)})$ (see the diagram below). As $\iota_1: \mathfrak S \ra \mathfrak S_{\log}^{(1)}$ is $\Gal_{K_\infty}$-equivariant where $\Gal_{K_\infty}$ acts trivially on $\mathfrak S$, we know that $\Gal_{K_\infty}$ acts trivially on $\fM$ as a submodule $\fM \hookrightarrow \fM \otimes_{\mathfrak S, \iota_1} \mathfrak S_{\log}^{(1)}$. Moreover, the action of $\tau$ comes from the descent data, in other words, we have a commutative diagram of isomorphisms 
\[
\begin{tikzcd}%[row sep = 1em]
\fM \otimes_{\mathfrak S, \iota_1} \mathfrak S_{\log}^{(1)} \arrow[r, "\alpha_1^{-1}"] \arrow[d, "\beta"] & \mE (\mathfrak S_{\log}^{(1)})  \arrow[d, "\tau"] \\
\fM \otimes_{\mathfrak S, \iota_2} \mathfrak S_{\log}^{(1)} \arrow[r, "\alpha_2^{-1}"] & \mE (\mathfrak S_{\log}^{(1)}) \arrow[r, "\alpha_1"] & 
\fM \otimes_{\mathfrak S, \iota_1} \mathfrak S_{\log}^{(1)}
\end{tikzcd}
\]
where $\beta$ is the descent isomorphism, and the action of $\tau$ is given by $\alpha_1 \circ \tau \circ \alpha_1^{-1} = \alpha_1 \circ \alpha_2^{-1} \circ \beta.$ In particular, for any pre-log prism $(A, I, M_A)$ with a map $\tau$-equivariant map $e: \mathfrak S_{\log}^{(1)} \ra A$ such that two compositions $e \circ p_1, e \circ p_2: \mathfrak S \ra A$ agree, the induced action of $\tau$ on $\mE (A)$ is trivial. Apply this to the map $\mathfrak S_{\log}^{(1)} \xrightarrow{\iota^{(1)}} \Ainf \ra W(\cl k)$, we know that the action $\tau - 1$ sends $\fM$ into $\fM \otimes_{\mathfrak S, \iota_1} I^+ \mathfrak S_{\log}^{(1)}$, therefore we obtain an object in $\tu{Mod}^{\phi, \hat G}_{/\mathfrak S}$ (resp. an object in $\tu{Mod}^{\phi, G_K}_{/\mathfrak S}$ by extending the $\Gal_K$-action to $\fM\otimes_{\mathfrak S, \iota_\varphi} \Ainf$ semi-linearly). To obtain the desired functor that takes values in $\tu{Mod}^{\phi, \tau_{\tu{st}}}_{/\mathfrak S}$ (resp. $\tu{Mod}^{\phi, \tau}_{/\mathfrak S}$), we need to show that the action of $\tau - 1$ on $\fM$ is ``divisible by $\mu$'' in $\fM_{\Acrys} := \fM \otimes_{\mathfrak S, \iota_\varphi} \Acrys$. To this end, we further evaluate $\mE$ on the prism $\mS$ (see Exampple \ref{example:Breuil}) and write $\mM = \mE (\mS) \cong \fM \otimes_{\mathfrak S, \varphi} \mS$. One may similarly construct the $\tau$-action on $\mM \otimes_{\mS} \Acrys \cong \fM \otimes_{\mathfrak S, \varphi} \Acrys$ which agrees with the action $\tau$ extended from above. Comparing with the construction of the monodromy operator on $\mM$, we have
\begin{equation} \label{eq:formula_of_tau_and_N}
(\tau - 1) \cdot x = \sum_{m \ge 1}  \mN^m (x) \otimes \frac{([\epsilon] - 1)^m}{m!} \in \mM \otimes_{\mS} \Acrys 
\end{equation}
for all $x \in \mM$ (also see \cite[5.1.2]{Tong_Liu_compositio}), which in particular lives in  $\fM \otimes_{\mathfrak S, \iota_{\varphi}} ([\epsilon] - 1) \Acrys$ (note that $[\epsilon -1]^p \in p \Acrys$).  

Next we show that both functors are equivalences. %It is in fact not difficult to show that these functors are fully faithful.  
To start, let us record the following variant of the category $\tu{Mod}^{\phi, N}_{\mathfrak S}$ considered in Subsection \ref{ss:BK_mod_with_int_monodromy}. Let $\breve K_0 = W(\cl k)[1/p]$. We denote by $\tu{Mod}^{\phi, N_{\breve K_0}}_{/\mathfrak S}$ the category of Breuil--Kisin modules $(\fM, \varphi) \in \tu{Mod}^{\phi}_{/\mathfrak S}$ equipped with a $\breve K_0$-linear monodromy operator $N$ on $\fM/u \fM \otimes_{W(k)}\breve K_0 $ satisfying the usual condition $N \varphi = p \varphi N$. Reverse engineering from  (\ref{eq:formula_of_tau_and_N}), we have a natural functor 
\[
D^{N_{\breve K_0 }}:  \tu{Mod}^{\phi, \tau}_{/\mathfrak S} \lra \tu{Mod}^{\phi, N_{\breve K_0}}_{/\mathfrak S}
\]
sending a $(\phi, \tau)$-module over $\mathfrak S$ to its underlying Breuil--Kisin module $(\fM, \varphi)$, equipped with a monodromy operator $N$ extended by linearity from 
\begin{equation} \label{eq:formula_of_N}
N (x) := \vartheta ( \frac{(\tau - 1) x}{[\epsilon] - 1}) \in \mM \otimes_{\mS} W(\cl k) \subset \fM/u\fM \otimes \breve K_0 
\end{equation}
for $x \in  \fM/u \fM$, where $\vartheta: \mM \otimes_{\mS} \Acrys \ra \mM \otimes_{\mS} W(\cl k)$ is the   map induced from   $\Acrys \ra W(\cl k)$. It is clear that $N \varphi = p \varphi N$. Moreover, the functor $D^{N_{\breve K_0}}$ is faithful and the composition with the functors in (\ref{eq:equivalence_for_phi_tau}) constructed above agrees with the composition 
\[
\tu{Rep}^{\tu{st}, \ge 0}_{\Z_p} (\tu{Gal}_K) \lra \tu{Mod}^{\phi, N_{\tu{int}}}_{/\mathfrak S} \lra  \tu{Mod}^{\phi, N_{\breve K_0}}_{/\mathfrak S}
\] 
which is fully faithful (the first functor is fully faithful by Theorem \ref{thm:D_Sigma_fully_faithful}). This shows that both functors in (\ref{eq:equivalence_for_phi_tau}) are fully faithful. To finish the proof, it suffices to show that each $(\phi, \tau)$-module in $ \tu{Mod}^{\phi, \tau}_{/\mathfrak S} $ comes from a $\Z_p$-lattice in a semistable $\Gal_K$-representation. Let $\fM \in \tu{Mod}^{\phi, \tau}_{/\mathfrak S}$. Note that we already have a $\Z_p$-valued $\Gal_K$ representation $T = (\fM \otimes W(C^\flat))^{\varphi = 1}$  coming from the ``\'etale realization'' of $\fM$. It suffices to show that $T[1/p]$ is semistable. Also note that, since $\Bst$-admissibility is insensitive to unramified base change of $K$, we may assume that $N$ from (\ref{eq:formula_of_N}) is already defined on $D:= \fM/u \fM \otimes \Q_p$. Now, the condition on $\fM$ ensures that the induced action of  $\tau$ on $D$ is trivial. Therefore, $(\fM \otimes \Bstp)^{\Gal_K}$ contains the following subspace 
\begin{equation} \label{eq:D_bar}
\cl D := \left\{ \sum_{i \ge 0} N^i (y) \otimes \frac{(\log [\varpi^\flat])^i}{i!} \in D \otimes \Bstp \cong \mM \otimes_{\mS} \Bstp  \: \vline \: y \in D  \right\}
\end{equation}
which has $K_0$-dimension equal to $\dim_{K_0} D$,\footnote{This observation is essentially due to T. Liu, see \cite[Section 7.2]{Tong_Liu_torsion}.} this shows that
\[
\dim_{K_0} (T \otimes_{\Z_p} \Bst)^{\Gal_K} = \dim_{K_0} ((\fM \otimes \Ainf[1/\mu]) \otimes_{\Ainf[1/\mu]} \Bst)^{\Gal_K} \ge \dim_{K_0} D
\]
from which we conclude that $T[1/p]$ is indeed semistable. The claim on the image of crystalline lattices is evident from the description of the monodromy operator in (\ref{eq:formula_of_N}).
\eproof

\br \label{remark:relation_to_other_phi_tau}
The categories $\tu{Mod}_{/\mathfrak S}^{\phi, \hat G}$, $\tu{Mod}_{/\mathfrak S}^{\phi, G_K}$, and a closely related category of $(\phi, \tau)$-modules over $\mathfrak S[\frac{1}{E}]^{\wedge}_p$ have already been studied in \cite{DuLiu, Gao_Hui, Caruso}.   %and is denoted by %$\tu{Mod}_{/\mathfrak S, \Ainf}^{\phi, G_K}$  therein, resp. in \cite{Caruso}). %Note that the conditions (1) and (2) in Definition \ref{definitions:other_similar_BK_modules} are phrased slightly differently in \cite{DuLiu, Gao_Hui}. 
%\bi \item $\fM \subset (\fM_{\Ainf})^{G_{\infty}}$ \item $\fM/u \fM \subset  $ {\color{red} This condition is confusing}\ei 
The main result of \cite{DuLiu} is that there is an equivalence between $ \tu{Rep}^{\tu{st}, \ge 0}_{\Z_p} $ and $\tu{Mod}^{\phi, \hat G}_{/\mathfrak S}$ (which they then use in their alternative proof of Theorem \ref{mainthm:equivalence}). Their proof of this equivalence uses an earlier theory $(\varphi, \hat G)$-modules over a rather mysterious ring $\widehat{\mathcal R} \subset \Ainf$ developed by T. Liu \cite{Tong_Liu_torsion}. On the other hand, the main result of \cite{Gao_Hui} states that there is an equivalence between $ \tu{Rep}^{\tu{st}, \ge 0}_{\Z_p} $ and $\tu{Mod}^{\phi, G_K}_{/\mathfrak S}$. The proof of this equivalence in \textit{loc.cit.} %is based on the theory of representations of finite $E$-height and
relies on the construction of the operator $N_{\nabla}$ (and showing that it is valued over a small enough ring). One construction of such a $N_{\nabla}$ appeared in \cite{Caruso} using $p$-adic approximation techniques but contained a serious gap. Gao \cite{Gao_Hui} bypassed this issue using some delicate analysis of locally analytic vectors. Our Theorem \ref{thm:equivalence_with_phi_tau} is a variant of these results in \cite{DuLiu, Gao_Hui} (in fact we expect that the categories $ \tu{Mod}^{\phi, \tau_{\tu{st}}}_{/\mathfrak S}$ and $\tu{Mod}_{/\mathfrak S}^{\phi, \tau}$ are more useful in applications to study  Galois representations, though \textit{a posteriori} that the categories in (\ref{diagram:several_BK_categories}) are all equivalent).  In contrast to the arguments from \cite{DuLiu, Gao_Hui, Caruso}, our proof above only uses the main equivalence (Theorem \ref{mainthm:equivalence}) and does not involve the ring $\widehat{\mathcal R}$, or the notion of locally analytic vectors, or any $p$-adic approximation. In particular, our proof does not fix the gap in \cite{Caruso}.  
\er 

\br[Essential image of $D_{\mathfrak S}$] \label{remark:essential_image_D_sigma}
The proof of Theorem \ref{thm:equivalence_with_phi_tau} in fact allows us to describe the essential image of $D_{\mathfrak S}$ in Theorem \ref{thm:D_Sigma_fully_faithful}. Given a $(\phi, N_{\tu{int}})$-module $\fM$ over $\mathfrak S$, we can \textit{construct} a $\tau$-action $\tau: \fM \ra \fM \otimes_{\mathfrak S, \varphi} \Acrys [1/p]$ by 
\[\tau (x) := \sum_{m \ge 0}  \mN^m (x) \otimes ( {([\epsilon] - 1)^m}/{m!}),\] where $\mN$ is given by $ N \otimes \tu{id} + \tu{id} \otimes \mN_{\mS}$ on $\fM$ as a submodule of $\fM \otimes_{\mathfrak S, \varphi} \mS [1/p] \cong (\fM/u \fM) \otimes_{W, \varphi} \mS [1/p]$ (also see \cite[Lemma 5.1.1]{Tong_Liu_compositio}). From the proof above, we know that $\fM = D_{\mathfrak S} (L)$ for some lattice $L \in \tu{Rep}^{\tu{st}}_{\Z_p} (\Gal_K)$ precisely when $\im (\tau) \subset \fM \otimes_{\mathfrak S} \Ainf$. The essential image of $\tu{Rep}^{\tu{crys}}_{\Z_p} (\Gal_K)$ in $\tu{Mod}^{\phi}_{/\mathfrak S}$ has a similar description, with $\mN = \tu{id}\otimes \mN_{\mS}$ in the formula above. 
\er 

%\br[Liu] The equivalences above can be obtained  \er 

\subsection{Lattices in filtered $(\phi, N)$-modules} We end this article by providing a simple construction of the functor $M_{\tu{st}}$ previously studied by T. Liu in \cite{TongLiu_integral_MF_lattice}. As in the introduction, we let $\tu{MF}^{\phi, N, \tu{int}}_{/K} $ denote the category of finite free $W(k)$-modules $M$ equiped with a $W(k)$-semilinear Frobenius operator $\varphi: M \ra M$ that becomes an isomorphism upon inverting $p$, a $W(k)$-linear monodromy operator $N: M \ra M$ satisfying $N \varphi = p \varphi N $ and a filtration $\Fil^\bullet D_K$ on $D_K = M \otimes_{W(k)} K.$ 

\bp \label{prop:TongLiu}
There exists a natural faithful functor 
\[M_{\tu{st}}: \tu{Rep}^{\tu{st}, \ge 0}_{\Z_p} (\tu{Gal}_K) \lra \tu{MF}^{\phi, N, \tu{int}}_{/K}
\]
giving rise to a $(\phi, N)$-stable lattice $M_{\tu{st}}(L)$ inside $D_{\tu{st}}(L[1/p])$ for every $L \in  \tu{Rep}^{\tu{st}, \ge 0}_{\Z_p} (\tu{Gal}_K).$ 
Moreover, when $er \le p-2$, $M_{\tu{st}}$ restricts to a fully faithful functor 
\[
\tu{Rep}^{\tu{st}, [0, r]}_{\Z_p} (\tu{Gal}_K) \lhook\joinrel\longrightarrow \tu{MF}^{\phi, N, \tu{int}}_{/K},
\]
where $ \tu{Rep}^{\tu{st}, [0, r]}_{\Z_p} (\tu{Gal}_K)$ denotes the subcategory of $\tu{Rep}^{\tu{st}}_{\Z_p} (\tu{Gal}_K)$
such that the Hodge--Tate weights of the corresponding representations are in $\{0, ..., r\}$, and $e$ denotes the ramification of $K/\Q_p$.
\ep 

\bproof 
The construction of $M_{\tu{st}}$ is immediate from either evaluting the corresponding prismatic $F$-crystal on the Hyodo--Kato log prism $W_{\tu{HK}}$ (see Example \ref{remark:Breuil_HK} and Remark \ref{remark:monodromy_on_D_0}), or composing the functor $D_{\mathfrak S}$ from Theorem \ref{thm:D_Sigma_fully_faithful} with the natural functor 
\begin{equation} \label{eq:functor_tilde_M_st}
\tu{Mod}^{\phi, N_{\tu{int}}}_{/\mathfrak S} \lra \tu{MF}^{\phi, N, \tu{int}}_{/K}
\end{equation}
which sends $(\fM, \varphi, N)$ to $(M, \varphi, N, \Fil^\bullet)$ where $M = \fM/u$ and $\Fil^\bullet M\otimes_{W} K$ is constructed using the Frobenius on $\fM$ as in the proof of Theorem \ref{thm:D_Sigma_fully_faithful}. The functor in (\ref{eq:functor_tilde_M_st}) is faithful by considering the induced functor on the isogeny categories. When $er \le p -2$, it is in fact full by an observation of T. Liu \cite[Proposition 2.16]{TongLiu_integral_MF_lattice}, using a  result on the structure of certain torsion $\mathfrak S$-modules (Lemma 2.18 in \textit{loc.cit.}). 
\eproof

\br 
The categories discussed in this section naturally fit into the following diagram, which also serves as a summary of this article.  
\begin{equation} \label{eq:the_summary_diagram} 
\begin{tikzcd}[column sep = 3.7em]%[ row sep = 2em]
& \tu{Mod}^{\phi,\hat G}_{/\mathfrak S} \arrow[r, "\sim"]& \tu{Mod}^{\phi,G_K}_{/\mathfrak S} \\
\tu{Vect}^{\varphi} (\mO_K^{\log})_{\Prism}^{\tu{eff}}   \arrow[r, "\sim"]  \arrow[r, swap, "\tu{Thm.}  \ref{thm:equivalence_with_phi_tau}"]  \arrow[d, swap, "\rotatebox{90}{$\sim$}"]  \arrow[d,  "\tu{Thm.} \ref{mainthm:equivalence}"]  & \tu{Mod}^{\phi, \tau_{\tu{st}}}_{/\mathfrak S} \arrow[r, "\sim"] \arrow[r, swap, "\tu{Thm.}  \ref{thm:equivalence_with_phi_tau}"] \arrow[u, "\rotatebox{-90}{$\sim$}"]  &   \tu{Mod}^{\phi, \tau}_{/\mathfrak S}  \arrow[d, hook]   \arrow[u, "\rotatebox{-90}{$\sim$}"]    \\
\tu{Rep}^{\tu{st}, \ge 0}_{\Z_p} (\tu{Gal}_K) %\arrow[rru, swap, near end, "\tu{Thm. }\ref{thm:equivalence_with_phi_tau}"] %\arrow[ur, "\rotatebox{40}{$\sim$}"] 
\arrow[rr, hook, "\tu{Thm. } \ref{thm:D_Sigma_fully_faithful}"]  
\arrow[rrd, near end, swap, "{\tu{Prop. } \ref{prop:TongLiu}}"] & & \tu{Mod}^{\phi, N_{\tu{int}}}_{/\mathfrak S} \arrow[r, hook] \arrow[d] 
& \tu{Mod}^{\phi, N}_{/\mathfrak S}. \\
& & \tu{MF}^{\phi, N, \tu{int}}_{/K} 
\end{tikzcd}
\end{equation}
Here in this diagram, we use $\tu{Vect}^{\varphi} (\mO_K^{\log})_{\Prism}^{\tu{eff}}$ to denote the \textit{effective} prismatic $F$-crystals on $(\mO_K, \varpi^{\N})$ for siplicity, namely $F$-crystals $\mE$ such that $\varphi_{\mE}(\varphi^* \mE) \subset \mE$, this is equivalent to $\tu{Rep}^{\tu{st}, \ge 0}_{\Z_p} (\tu{Gal}_K)$ by Remark \ref{remark:HT_weights}.  
\er 
%$\tu{Mod}^{\phi, N}_{/\mathfrak S}$ denote the category of Breuil--Kisin modules $(\fM, \varphi)$ equipped with a monodromy operator $N$ on $\fM/u \otimes \Q_p$ introduced in \cite{Kisin_crystal}. %All other categories and arrows will be introduced in Section \ref{section:BK}
%we use $\hookrightarrow$ to denote fully faithful embeddings 

%\bibliographystyle{alpha}
\bibliographystyle{plain}
\bibliography{ref}
\end{document}